\documentclass[12pt,twoside]{amsart}
\usepackage[margin=3cm]{geometry}
\usepackage[colorlinks=false]{hyperref}
\usepackage[english]{babel}
\usepackage{graphicx,titling}
\usepackage{float}
\usepackage{amsmath,amsfonts,amssymb,amsthm}
\usepackage{lipsum}
\usepackage[T1]{fontenc}
\usepackage{fourier}
\usepackage{color}
\usepackage[latin1]{inputenc}
\usepackage{esint}
\usepackage{caption}
\usepackage{ccicons}

\makeatletter
\def\blfootnote{\xdef\@thefnmark{}\@footnotetext}
\makeatother

\newcommand\ccnote{
    \blfootnote{\ccLogo\, \ccAttribution\,\, Licensed under a Creative Commons Attribution License (CC-BY).}
}

\usepackage[export]{adjustbox}
\numberwithin{equation}{section}
\usepackage{setspace}
\renewcommand{\le}{\leqslant}

\renewcommand{\ge}{\geqslant}

\renewcommand{\mathbb}{\varmathbb}
\usepackage{fancyhdr}
\def\R{\mathbb{R}}     \def\Sph{\mathbb{S}} 
\def\C{\mathbb{C}}   
\renewcommand{\epsilon}{\varepsilon}
\def\sing{{\rm sing\,\,}}
\def\dvg{{\rm div}}
\def\dist{{\rm dist\,}}
\def\graph{{\rm graph\,}}
\def\spt{{\rm spt\,}}
\def\res{\hbox{ {\vrule height .22cm}{\leaders\hrule\hskip.2cm} } }        
\def\sres{\hskip-0.7pt\raise-0.5pt\hbox{ {\vrule height .18cm}{\leaders\hrule\hskip.17cm} }\hskip-0.7pt }        
\def\wtilde#1{\,\,\widetilde{\hskip-2pt M}}
\def\thn{\hskip1.2pt}

\def\tsum{{\textstyle\sum}}
\def\tint{{\textstyle\int}}
\def\fr#1#2{{\textstyle\frac{#1}{#2}}}
\def\smfr#1#2{{\scriptsize\frac{#1}{#2}}}

\def\fr#1#2{{\textstyle\frac{#1}{#2}}}
\def\ha{\fr{1}{2}} 
\def\avint{\hskip1pt\text{\small $-$}\hskip-10.6pt\int}                                              
\def\smavint{\text{\scriptsize $-$}\hskip-6.85pt\tint}
\DeclareSymbolFont{upright}       {OMS}{zplm}{m}{n}
\SetSymbolFont{upright}{bold}     {OMS}{zplm}{b}{n}
\DeclareMathSymbol{\emptyset}{\mathord}{upright}{"3B}  

\def\tg#1{\tag*{\textrm{#1}}} 
\newcounter{sequation}[section] 
\renewcommand{\thesequation}{\arabic{section}.\arabic{sequation}} 
 \def\tl#1{\refstepcounter{sequation}\label{#1}{\textbf{\thesequation}\thn}} 
 \def\dl#1{\refstepcounter{sequation}\label{#1}\leqno\textrm{\thesequation}} 
\def\dtg#1{\refstepcounter{sequation}\label{#1}\tag*{\textrm{\thesequation}}} 
\newcounter{pequation} 
\renewcommand{\thepequation}{(\arabic{pequation})}%
\def\pdl#1{\refstepcounter{pequation}\label{#1}\leqno\textrm{\thepequation}} 
\def\ptg#1{\refstepcounter{pequation}\label{#1}\tag*{\textrm{\thepequation}}}     

\newenvironment{state}[1]{\par\smallbreak\noindent {\bf#1} %
\textit\bgroup}{\egroup \par\ifdim\lastskip<\medskipamount \removelastskip\penalty55\medskip\fi} %
\renewenvironment{proof}{\setcounter{pequation}{0}\noindent}{\,\,\nolinebreak$\Box$} 

\address{Leon Simon\\ Mathematics Department\\ Stanford University\\ Stanford CA 94305}
\email{lsimon@stanford.edu}

\begin{document}

\thispagestyle{empty}

\begin{minipage}{0.28\textwidth}
\begin{figure}[H]
\includegraphics[width=2.5cm,height=2.5cm,left]{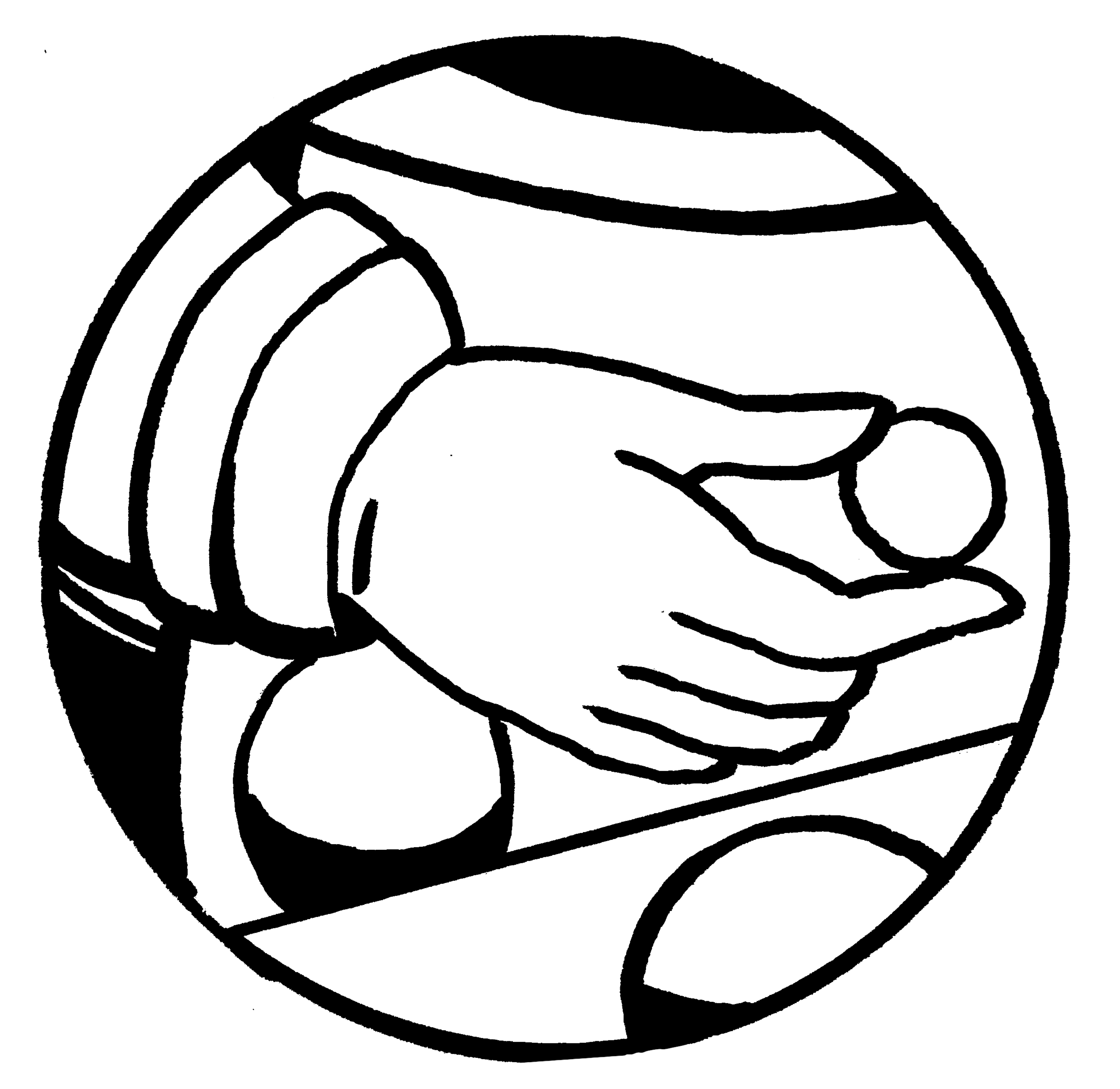}
\end{figure}
\end{minipage}
\begin{minipage}{0.7\textwidth} 
\begin{flushright}
Ars Inveniendi Analytica (2021), Paper No. 5, 35 pp.
\\
DOI 10.15781/hk9g-zz18
\end{flushright}
\end{minipage}

\ccnote

\vspace{1cm}


\begin{center}
\begin{huge}
\textit{        A Liouville-type theorem for\\  stable minimal  hypersurfaces }
\end{huge}
\end{center}

\vspace{1cm}


\begin{center}
\large{\bf{Leon Simon}} \\
\vskip0.15cm
\footnotesize{Stanford University}
\end{center}

\vspace{1cm}


\begin{center}
\noindent \em{Communicated by Guido De Philippis}
\end{center}
\vspace{1cm}


\noindent \textbf{Abstract.} \textit{We prove that if $M$ is a strictly stable complete minimal hypersurface in
$\smash{\R^{n+1+\ell}}$ which has finite density at infinity and which lies on one side of a cylinder
$\smash{\C=\C_{0}\times\R^{\ell}}$, where $\smash{\C_{0}}$ is a strictly stable area minimizing hypercone in
$\smash{\R^{n+1}}$ with $\sing\C_{0}=\{0\}$, then $M$ must be cylindrical---i.e.\ $M=S\times\R^{\ell}$, where
$S\subset\smash{\R^{n+1}}$ is a smooth strictly stable minimal hypersurface in $\,\smash{\R^{n+1}}$\!.  \,Applications will be
given in~\emph{\cite{Sim21a},~\cite{Sim21b}}. } \vskip0.3cm

\vspace{0.5cm}



\section*{Introduction}\label{intro}

\noindent The main theorem here (Theorem~\ref{main-th}) establishes that if $M$ is a complete strictly stable minimal
hypersurface in $\R^{n+1+\ell}$ lying on one side of a cylindrical hypercone $\C=\C_{0}\times\R^{\ell}$ with
$\sing\C_{0}=\{0\}$ and $\C_{0}$ strictly stable and minimizing, then $M$ must itself be cylindrical---i.e.\ of the form
$S\times\R^{\ell}$ where $S$ is a smooth complete (not necessarily connected) strictly stable minimal hypersurface in $\R^{n+1}$.

This result, or more correctly its Corollary~\ref{co-1} below, is a crucial ingredient in the author's recent proof
(in~\cite{Sim21b}) that, with respect to a suitable $C^{\infty}$ metric for $\R^{n+1+\ell}$, $n+\ell\ge 8$, there are
examples of strictly stable minimal hypersurfaces which have singular set of the form $\{0\}\times K$, where $K$ is an
arbitrary closed subset of $\R^{\ell}$.

An outline of the present paper is as follows. After a description of the main results in~\S1 and some notational and technical
preliminaries in \S2 and \S3, in \S4 we establish $L^{2}$ growth estimates for solutions of the Jacobi equation (i.e.\ the
linearization of the minimal surface equation) on $M$. These estimates are applied to give growth estimates on (i) $(x,y)\cdot\nu$,
where $\nu$ is the unit normal of $M$, (ii) $\nu_{\!y}=(e_{n+2}\cdot\nu,\ldots,e_{n+1+\ell}\cdot\nu)$ (i.e.\ the components of
the unit normal $\nu$ of $M$ in the $y$-coordinate directions), and  (iii) (in \S5) $d|M$, where $d(x)=\dist((x,y),\C)$ is the
distance to the cylinder $\C$. In~\S6 the growth estimates established in \S4 and \S5 are combined to show that, if $\nu_{\!y}$ is
not identically zero,
\[%
R^{\gamma-\alpha}\le \int_{M\cap \{(x,y):|(x,y)|<R\}}\nu_{\!y}^{2}\,d\mu \le  %
R^{-2+\gamma+\alpha}, \,\, \gamma=\ell+2+\beta_{1} ,\,\,  %
\]%
for each $\alpha\in (0,1)$ and all sufficiently large $R$ (depending on $\alpha$ and $M$), where $\mu$ is $(n+\ell)$-dimensional
Hausdorff measure, and $\beta_{1}=2((\fr{n-2}{2})^{2}+\lambda_{1})^{1/2}$, with $\lambda_{1}$ the first eigenvalue of the
Jacobi operator of the compact minimal surface $\Sigma=\C_{0}\cap \Sph^{n}$.  These inequalities are clearly impossible for
$R>1$, so we finally conclude that indeed $\nu_{\!y}$ is identically zero, showing that $M$ is cylindrical as
claimed in the main theorem.

\section{Main Results}\label{main-res}

\noindent Let $\C_{0}$ be a minimal hypercone in $\R^{n+1}$ with $\sing\C_{0}=\overline\C_{0}\setminus\C_{0}=\{0\}$
and let $\mathcal{L}_{\C_{0}}$ be the Jacobi operator (linearized minimal surface operator) defined by
\[%
\mathcal{L}_{\C_{0}}u=\Delta_{\C_{0}}u+|A_{\C_{0}}|^{2}u,
\]%
where $|A_{\C_{0}}|^{2}$ is the squared length of the second fundamental form of $\C_{0}$. In terms of the Jacobi operator
\[%
\mathcal{L}_{\Sigma}=\Delta_{\Sigma}u+|A_{\Sigma}|^{2}u
\dl{jac-op-Sig}
\]%
of the compact submanifold $\Sigma=\C_{0}\cap \Sph^{n}$ (which is minimal in $\Sph^{n}$), we have
\[%
\mathcal{L}_{\C_{0}}u=r^{1-n}\frac{\partial}{\partial r}\bigl(r^{n-1}\frac{\partial u}{\partial r}\bigr) %
+r^{-2}\mathcal{L}_{\Sigma},\quad r=|x|.  %
\dl{jac-op}
\]%
If $\lambda_{1}$ is the first eigenvalue of $-\mathcal{L}_{\Sigma}$, and $\varphi_{1}>0$ is the corresponding eigenfunction
(unique up to a constant factor) we have (possibly complex-valued) solutions $r^{\gamma_{1}^{\pm}}\varphi_{1}$ of
$\mathcal{L}_{\C_{0}}u=0$, where $\gamma_{1}^{\pm}$ are the ``characteristic exponents''
\[%
\gamma_{1}^{\pm} = -\frac{n-2}{2} \pm \sqrt{ \bigl(\frac{n-2}{2}\bigr)^{2}+\lambda_{1}}.  %
\dl{ch-exps}
\]%
We assume $\C_{0}$ is strictly stable, i.e.\ there is $\lambda>0$ with
\[%
\lambda\int_{\C_{0}}|x|^{-2}\zeta^{2}(x)\,d\mu(x)\le %
\int_{\C_{0}}\bigl( |\nabla_{\C_{0}}\zeta|^{2} - |A_{\C_{0}}|^{2}\zeta^{2}\bigr)\,d\mu\,\,\,\, %
\forall \zeta\in C_{c}^{\infty}(\R^{n+1}),%
\dl{st-stab}
\]%
where $\mu$ is $n$-dimensional Hausdorff measure; subsequently $\mu$ will always denote Hausdorff measure of the
appropriate dimension.

Using~\ref{jac-op} and \ref{ch-exps}, \ref{st-stab} is readily checked to be equivalent to the condition that
\[%
\lambda_{1}>-((n-2)/2)^{2}, %
\dl{st-stab-1}
\]%
in which case
\[%
\gamma_{1}^{-}<-\frac{n-2}{2} <\gamma_{1}^{+}<0.
\dl{ch-ineq}
\]%
The main theorem here relates to hypersurfaces  $M\subset\R^{n+1+\ell}$, where $\ell\ge 1$; points in 
$\R^{n+1+\ell}=\R^{n+1}\times\R^{\ell}$ will be denoted
\[%
(x,y)=(x_{1},\ldots,x_{n+1},y_{1},\ldots,y_{\ell}). 
\]%

In the main theorem, which we now state and which will be proved in~\S\ref{th-1-pf}, we assume that~(i) $M$ is a smooth
complete minimal hypersurface in $\R^{n+1+\ell}$ lying on one side of $\C=\C_{0}\times\R^{\ell}$, i.e.\
\[%
M\subset U_{+}\times\R^{\ell}, %
\dl{one-side}
\]%
where $U_{+}$ is one of the two connected components of $\R^{n+1+\ell}\setminus\overline{C}_{0}$, that (ii) $M$ is strictly
stable in the sense that (Cf.\ \ref{st-stab}) there is $\lambda>0$ with
\[%
\lambda\int_{M}|x|^{-2}\zeta^{2}(x,y)\,d\mu(x,y)\le %
\int_{M}\bigl( |\nabla_{M}\zeta|^{2} - |A_{M}|^{2}\zeta^{2}\bigr)\,d\mu\,\,\,\, %
\dl{str-stab-M} %
 \]%
 for all $\zeta\in C_{c}^{1}(\R^{n+1+\ell})$, and that (iii) $M$ has finite density at $\infty$, i.e.\
\[%
\sup_{R>1}R^{-n-\ell}\mu(M\cap B_{R})<\infty,  %
\dl{fin-dens}
\]%
where, here and subsequently, $B_{R}$ is the closed ball of radius $R$ and centre $(0,0)$ in $\R^{n+1+\ell}$:
\[%
B_{R}=\{(x,y)\in\R^{n+1+\ell}:|(x,y)|\le R\};
\]%
the corresponding open ball, which we shall occasionally use, is denoted
\[%
\breve B_{R}=\{(x,y)\in\R^{n+1+\ell}:|(x,y)|< R\}.
\]%

\begin{state}{\bf{}\tl{main-th} Theorem (Liouville-type Theorem.)}%
If $M$ is a smooth, complete, embedded minimal hypersurface (without boundary) in \,$\R^{n+1+\ell}$ such
that~{\rm\ref{one-side}}, {\rm\ref{str-stab-M}} and  {\rm\ref{fin-dens}} hold, then $M$ is cylindrical, i.e.\
\[%
M=S\times\R^{\ell}, 
\]%
where $S$ is a smooth complete strictly stable minimal hypersurface which is contained in~$U_{+}$.
\end{state}%

\noindent{\bf{}Remark: } 
Using the regularity theory of~{\rm\cite{SchS81}}, the conclusion continues to hold, with no essential change in the proof, in case
$M$ is allowed to have a singular set of finite $(n+\ell-2)$-dimensional Hausdorff measure. Indeed if we use the regularity theory
of~{\rm\cite{Wic14}} then we need only assume \emph{a priori} that $M$ has zero $(n+\ell-1)$-dimensional Hausdorff measure.

\smallskip

For applications in~\cite{Sim21a}, \cite{Sim21b} we now state a corollary of the above theorem, in which the stability
hypothesis on $M$ is dropped, and instead $\C_{0}$ is assumed to be both strictly stable and strictly minimizing, and we
impose an \emph{a priori} smallness assumption on the last $\ell$ components of the unit normal
$\nu_{M}=(\nu_{1},\ldots,\nu_{n+1+\ell})$ of $M$; i.e.\ a smallness assumption on $(\nu_{n+2},\ldots,\nu_{n+1+\ell})\,$
which we subsequently write as
\[%
\nu_{y}=(\nu_{y_{1}},\ldots,\nu_{y_{\ell}});\,\, \text{  i.e.\   }  \nu_{y_{j}}=e_{n+1+j}\cdot\nu_{M},\,\,\, j=1,\ldots,\ell.
\dl{nu-y}
\]%
Notice that such a smallness assumption on $|\nu_{y}|$ amounts to a restriction on how rapidly $M$ can vary in the
$y$-directions.

Recall (see \cite{HarS85}) $\C_{0}$ is said to be strictly minimizing if there is a constant $c>0$ such that
{\abovedisplayskip8pt\belowdisplayskip8pt%
\begin{align*}%
  &\mu(\C_{0}\cap B_{1}) \le \mu(T\cap B_{1}) -c\rho^{n} \text{ whenever $\rho\in(0,\ha]$ and $T$ is a smooth
  }\dtg{st-min}\\
  \noalign{\nobreak\vskip-3pt}  %
&\hskip0.5in \text{compact hypersurface-with-boundary in $B_{1}\setminus B_{\rho}$ with  %
                                                                                                      $\partial T=\C_{0}\cap \Sph^{n}$}.  %
\end{align*}}%

\begin{state}{\bf{}\tl{co-1} Corollary.}%
Suppose $\C_{0}$ is strictly stable and strictly minimizing (as in~{\rm\ref{st-min}}) and $\alpha\in (0,1)$. Then there is
$\epsilon_{0}=\epsilon_{0}(\C_{0},\alpha)\in (0,\ha]$ such that if $M$ is a smooth, complete, embedded minimal
hypersurface in $\R^{n+1+\ell}$ with
\[%
\left\{\begin{aligned}%
&M\subset U_{+}\times\R^{\ell}, \\  %
\noalign{\vskip-2pt} %
&{\sup}_{R>1}R^{-n-\ell}\mu(M\cap B_{R})\le (2-\alpha)\mu(\C\cap B_{1}), and \\ %
\noalign{\vskip-2pt} %
&{\sup}_{M}|\nu_{y}| < \epsilon_{0}, 
\end{aligned}\right.%
\]%
then 
{\abovedisplayskip-3pt\belowdisplayskip8pt%
\[%
M=\lambda S\times\R^{\ell} 
\]}%
for some $\lambda>0$, where $S$ is the minimal hypersurface in $U_{+}$ as in~\emph{\cite{HarS85}}
(see~{\rm\ref{props-S}} in~{\rm\S\ref{proof-liou}} below).
\end{state}%

Applications of the above results are given in~\cite{Sim21a} and \cite{Sim21b}.  Although the assumptions of strict stability
in Theorem~\ref{main-th} and $|\nu_{y}|$ small in Corollary~\ref{co-1} are appropriate for the applications
in~\cite{Sim21a}, \cite{Sim21b}, it would be of interest to know if these restrictions can be significantly relaxed---for
example the question of whether or not mere stability would suffice in place of the strict stability assumption in
Theorem~\ref{main-th}.

\section[Preliminaries concerning $M$]{Preliminaries concerning $M$}\label{prelims}

\noindent As in~\S\ref{main-res}, $\C_{0}$ will be a smooth embedded minimal hypercone in $\R^{n+1}$ with
$\sing\C_{0}=\{0\}$, and we let $U_{\pm}$ be the two connected components of $\R^{n+1}\setminus \overline{\C}_{0}$.

$M$ will be a smooth embedded minimal hypersurface, usually contained in $U_{+}\times \R^{\ell}$, although some
results are formulated to apply locally in a ball, and also independent of the inclusion assumption $M\subset
U_{+}\times\R^{\ell}$.

For the moment we assume only that $M$ is minimal (i.e.\ has first variation zero) in $\breve B_{R}$. Thus
{\abovedisplayskip8pt\belowdisplayskip8pt%
\[%
\int_{M}\dvg_{M}Z(x,y)\,d\mu(x,y) =0, \,\, Z=(Z_{1},\ldots,Z_{n+1+\ell}),\,Z_{j}\in C^{1}_{c}(\breve B_{R}), %
\dl{stationarity}
\]%
where $\dvg_{M}Z$ is the ``tangential divergence'' of $Z$:
\begin{align*}%
\dvg_{M}Z_{|(x,y)}=\tsum_{j=1}^{n+\ell}\tau_{j}\cdot D_{\tau_{j}}Z %
&= \tsum_{k,m=1}^{n+1+\ell}\tsum_{j=1}^{n+\ell}\tau_{j}^{k}\tau_{j}^{m}D_{k}Z_{m} \\ %
&=\tsum_{i,j=1}^{n+1+\ell}g^{ij}(x,y)D_{i}Z_{j}(x,y),
\end{align*}}%
with $\tau_{j}=(\tau_{j}^{1},\ldots,\tau_{j}^{n+1+\ell})$, $j=1,\ldots,n+\ell$, any locally defined orthonormal basis of
$T_{(x,y)}M$, $(x,y)\in M$ and
\[%
g^{ij} = \delta_{ij}-\nu_{i}\nu_{j},\quad \nu=(\nu_{1},\ldots,\nu_{n+1+\ell}) \text{ a unit normal for } M.  %
\dl{g-ij}
\]%
For $v\in C^{2}(M)$ we let $\graph v$ be the graph of $v$ taken off $M$:
\[%
\graph v =\bigl\{(x,y)+v(x,y)\nu_{\!M}(x,y):(x,y)\in M\bigr\}
\]%
(notice that this may fail to be an embedded hypersurface unless $v$ has small enough $C^{2}$ norm), and we take
\[%
\mathcal{M}_{\!M}(v)= \text{ the mean curvature operator on $M$}. %
\dl{def-script-M}
\]%
Thus $\mathcal{M}_{\!M}(v)$ is the Euler-Lagrange operator of the area functional $\mathcal{A}_{\!M}$ on $M$,
defined by
\[%
\mathcal{A}_{\!M}(v) = \int_{M}J_{M}(V)\,d\mu
\]%
where $V(x,y) =(x,y)+v(x,y)\nu_{M}(x,y)$ is the graph map taking $M$ to $\graph v$ and $J_{M}(V)$ is the Jacobian
of $V$:
\[%
J_{M}(V)=\sqrt{\det\bigl(D_{\tau_{i}}V\cdot %
  D_{\tau_{j}}V\bigr)}=\sqrt{\det\bigl(\delta_{ij}+v_{i}v_{j}+v^{2}\tsum_{k}h_{ik}h_{jk}   +2vh_{ij}\bigr)}, %
\]%
where, for $(x,y)\in M$, $\tau_{1},\ldots,\tau_{n+\ell}$ is an orthonormal basis for $T_{(x,y)}M$, $h_{ij}$ is the second
fundamental form of $M$ with respect to this basis, and $v_{ i}=D_{\tau_{i}}v$. Since $\sum_{i=1}^{n+\ell}h_{ii}=0$ we then
have
\[%
J_{M}(V) =  \sqrt{1+|\nabla_{M}v|^{2}-|A_{M}|^{2}v^{2}+E\bigl(v(h_{ij}),v^{2}(\tsum_{k}h_{ik}h_{jk}),(v_{i}v_{j})\bigr)}
\dl{jac-M}
\]%
where $E$ is a polynomial of degree $n+\ell$ in the indicated variables with constant coefficients depending only on $n,\ell$
and with each non-zero term having at least degree $3$.  So the second variation of $\mathcal{A}_{M}$ is given by
\[%
\fr{d^{2}}{dt^{2}}\bigl|_{t=0}\mathcal{A}(t\zeta)=\int_{M}\bigl(|\nabla_{M}\zeta|^{2}-|A_{M}|^{2}\zeta^{2}\bigr)\,d\mu%
=-\int_{M}\zeta\mathcal{L}_{M}\zeta\,d\mu, \quad \zeta\in C^{1}_{c}(\R^{n+1+\ell}),
\]%
where $\mathcal{L}_{M}$ is the  Jacobi operator on $M$ defined by
$$
\mathcal{L}_{M}v  = \Delta_{M}v+|A_{M}|^{2}v, %
\dl{L-M}
$$
with $\Delta_{M}f=\dvg_{M}(\nabla_{M}f)$ the Laplace-Beltrami operator on $M$ and $|A_{M}|=\bigl(\sum
h_{ij}^{2}\bigr)^{1/2}$ (the length of the second fundamental form of $M$). Of course, by~\ref{def-script-M} and
\ref{jac-M}, $\mathcal{L}_{M}$ is the linearization of the mean curvature operator at $0$, i.e.\
{\abovedisplayskip8pt\belowdisplayskip8pt%
\[%
\mathcal{L}_{M}v=\fr{d}{dt}\bigl|_{t=0}\mathcal{M}_{M}(tv).
\dl{L-M-2}
\]%
In view of~\ref{jac-M}  and the definition of $\mathcal{L}_{M}$, if $|A_{M}||u|+|\nabla_{M}u|<1$   on a domain
in $M$ then we can write  
\[%
\mathcal{M}(u)=\mathcal{L}_{M}(u) + \dvg_{M}E(u) + F(u) %
\dl{L-M-3}
\]}%
on that domain, where $|E(u)|\le C(|A_{M}|^{2} |u|^{2}+|\nabla_{M} u|^{2})$ and $|F(u)|\le
C(|A_{M}|^{3}v^{2}+|A_{M}||\nabla_{M}v|^{2})$, $C=C(n,\ell)$.

We say that $M$ is \emph{strictly stable} in  $\breve B_{R}$ if the second variation of $\mathcal{A}_{\!M}$ is strictly
positive in the sense that
\[%
{\inf}_{\zeta\in C^{1}_{c}(\breve B_{R}),\,\int_{M}|x|^{-2}\zeta^{2}(x,y)\,d\mu(x,y)=1} %
\,\,\,\fr{d^{2}}{dt^{2}}\bigl|_{t=0}\mathcal{A}(t\zeta) > 0, %
\]%
or, equivalently, that there is $\lambda>0$ such that~\ref{str-stab-M} holds.

Observe that if we have such strict stability then by replacing $\zeta$ in~\ref{str-stab-M} with $\zeta w$ we obtain
\begin{align*}%
&\smash[b]{\lambda\int_{M}|x|^{-2}w^{2}\zeta^{2}(x,y)\,d\mu(x,y)\le %
                \int_{M}\bigl( w^{2}|\nabla_{M}\zeta|^{2}+\zeta^{2}|\nabla_{M}w|^{2}}\\  %
\noalign{\vskip-1pt} %
&\hskip2.5in             +2\zeta w\nabla_{M}\zeta\cdot\nabla_{M}w -|A_{M}|^{2}\zeta^{2}w^{2}\bigr)\,d\mu, %
\end{align*}%
and $\int_{M}2\zeta w\nabla_{M}\zeta\cdot\nabla_{M}w=\int w\nabla_{M}\zeta^{2} \cdot\nabla_{M}w= -\int_{M}
(|\nabla_{M}w|^{2}\zeta^{2}+\zeta^{2}w\Delta_{M}w)$, so if $w$ is a smooth solution of $\mathcal{L}_{M}w=0$ on
$M\cap \breve B_{R}$ then
\[%
\lambda\int_{M}|x|^{-2}w^{2}\zeta^{2}\,d\mu(x,y)\le  \int_{M} w^{2}|\nabla_{M}\zeta|^{2}\,d\mu, %
\quad \zeta\in C^{1}_{c}(\breve B_{R}).  %
\dl{w-1-2-est} 
\]%
We need one further preliminary for $M$, concerning asymptotics of $M$ at $\infty$ in the case when $M$ is complete in
all of $\R^{n+1+\ell}$ with $M\subset U_{+}\times\R^{\ell}$:
\begin{state}{\bf{}\tl{tangent-cone} Lemma.}%
If $M$ is a complete embedded minimal hypersurface in all of \,$\R^{n+1+\ell}$ satisfying {\rm\ref{one-side}} and
{\rm\ref{fin-dens}}, then $M$ has $\C$ with some constant integer multiplicity $q$ as its unique tangent cone at $\infty$.

Furthermore if $M$ is stable (i.e.\ \emph{\ref{str-stab-M}} holds with $\lambda=0$), then for each $\delta>0$ there is
$R_{0}=R_{0}(\C_{0},q,\delta)>1$ such that
\[%
M\setminus \bigl(B_{R_{0}}\cup\bigl\{(x,y):|x|\le \delta|y|\bigr\}\bigr)\subset \cup_{j=1}^{q}\graph u_{j} \subset M 
\leqno{\rm(i)}
\]%
where each $u_{j}$ is a $C^{2}$ function on a domain $\Omega$ of \,$\C$ containing $\C\setminus
(B_{R_{0}}\cup\bigl\{(x,y):|x|\le \delta|y|\bigr\}$ and $\graph u_{j}=\{x + u_{j}(x,y)\nu_{\C_{0}}(x) :(x,y)\in\Omega\}$
($\nu_{\C_{0}}=$ unit normal of \,$\C_{0}$ pointing into $U_{+}$), and
\[%
\text{$\lim_{R\to\infty}{\sup}_{(x,y)\in \Omega\setminus
B_{R}}\bigl(|x||\nabla_{\C}^{2}u_{j} (x,y)|+|\nabla_{\C}u_{j}(x,y)|+|x|^{-1}u_{j}(x,y)\bigr)=0$.} %
\leqno{\rm(ii)}%
\]%
\end{state}%

\smallskip

\begin{proof}{\bf{}Proof:} Let $C(M)$ be a tangent cone of $M$ at $\infty$. Thus $C(M)$ is a stationary integer
multiplicity varifold with $\lambda C(M)=C(M)$ for each $\lambda>0$, and support of $C(M)\subset
\overline U_{+}\times\R^{\ell}$ and, by the compactness theorem of~\cite{SchS81}, $C(M)$ is stable and $\sing C(M)$
has Hausdorff dimension $\le n+\ell-7$. So by the maximum principle of~\cite{Ilm96} and the constancy theorem,
$C(M)=\C$ with constant multiplicity $q$ for some $q\in\{1,2,\ldots\}$.

Thus $\C$, with multiplicity $q$, is the unique tangent cone of $M$ at $\infty$, and the ``sheeting
theorem''~\cite[Theorem 1]{SchS81} is applicable, giving~(i) and~(ii).
\end{proof}

\section[Preliminaries concerning $\C_{0}$ and $\C$]{Preliminaries concerning $\C_{0}$ 
and $\C$}\label{prelim-C}

\noindent $\C_{0}\subset\R^{n+1}\setminus\{0\}$ continues to denote a smooth embedded minimal hypercone with
$\sing\hskip0.5pt \C_{0} =\overline\C_{0}\setminus\C_{0}=\{0\}$, $U_{\pm}$ denote the two connected components of
$\R^{n+1}\setminus\overline\C_{0}$, and we assume here that $\C_{0}$ is strictly stable as in~\ref{st-stab}.

With $\mathcal{L}_{\Sigma}$ the Jacobi operator of $\Sigma=\C_{0}\cap \Sph^{n}$ as in~\ref{jac-op-Sig}, we let
$\varphi_{1}>0,\,\varphi_{2},\ldots$ be a complete orthonormal set of eigenfunctions of $-L_{\Sigma}$. Thus
\begin{align*}%
& -L_{\Sigma}(\varphi_{j}) = \lambda_{j}\varphi_{j},\,\,\, j=1,2,\ldots, \text{with $\varphi_{1}>0$ and} \dtg{eigenvals}\\  %
\noalign{\vskip-1pt} %
&\hskip1.1in  \lambda_{1}<\lambda_{2}\le \lambda_{3}\le\cdots \le \lambda_{j}\le \cdots,\,\,
\lambda_{j}\to\infty\text{ as }j\to \infty, %
\end{align*}%
and every $L^{2}(\Sigma)$ function $v$ can be written $v=\sum_{j}\langle v,\varphi_{j} \rangle_{L^{2}(\Sigma)}\varphi_{j}$.

Notice that by orthogonality each $\varphi_{j},\,j\neq 1$, must then change sign in $\Sigma$, and, with
$\mathcal{L}_{\C_{0}}$ the Jacobi operator for $\C_{0}$ as in~\ref{jac-op},
\[%
\mathcal{L}_{\C_{0}}(r^{\gamma}\varphi_{j})=r^{\gamma-2}(\gamma^{2}+(n-2)\gamma-\lambda_{j})\,\varphi_{j}, %
\dl{exp-gamma}
\]%
so in particular if $\C_{0}$ is strictly stable (i.e.\ if $\lambda_{1}>-(\fr{n-2}{2})^{2}$ as in~\ref{st-stab-1}) we have
\begin{align*}%
&\mathcal{L}_{\C_{0}}(r^{\gamma_{j}^{\pm}}\varphi_{j}) =0, \quad \gamma_{j}^{\pm} = -\fr{n-2}{2} %
\pm\bigl((\fr{n-2}{2})^{2}+\lambda_{j}\bigr)^{1/2}, \dtg{char-exps}\\  %
\noalign{\vskip-2pt} %
&\hskip20pt -\infty\leftarrow\gamma_{j}^{-} \le \cdots\le \gamma_{2}^{-} <  %
 \gamma_{1}^{-}<-\fr{n-2}{2}<\gamma_{1}^{+}<\gamma_{2}^{+}\le \cdots\le  \gamma_{j}^{+}\to \infty.  %
\end{align*}%
Henceforth we write
\[%
\gamma_{j} =\gamma_{j}^{+},\quad j=1,2,\ldots. %
\dl{gamma-j}
\]%
The Jacobi operator $\mathcal{L}_{\C}$ on the minimal cylinder $\C=\C_{0}\times\R^{\ell}$ is
\[%
\mathcal{L}_{\C}(v) =  \mathcal{L}_{\C_{0}}(v)+\tsum_{j=1}^{\ell}D_{y_{j}}^{2}v,  %
\dl{jac-C}
\]%
and we can decompose solutions $v$ of the equation $\mathcal{L}_{\C}v=0$ in terms of the eigenfunctions $\varphi_{j}$
of~\ref{eigenvals}: for $v$ any given smooth function on $\C\cap \breve B_{R}$, we can write
$v(r\omega,y)=\sum_{j}v_{j}(r,y)\varphi_{j}(\omega)$, where
\[%
v_{j}(r,y)=\int_{\Sigma}v(r\omega,y)\,\varphi_{j}(\omega)\,d\mu(\omega), %
\dl{def-v-j}
\]%
and then $\mathcal{L}_{\C}v=0$ on $\breve B_{R}$ if and only if $v_{j}$ satisfies
\[%
r^{1-n}\frac{\partial}{\partial r}\bigl(r^{n-1}\frac{\partial v_{j}}{\partial r}\bigr) + %
\tsum_{k=1}^{\ell}\frac{\partial^{2}v_{j}}{\partial y_{k}^{2}} -\frac{\lambda_{j}}{r^{2}}v_{j}=0
\]%
for $(r,y)\in \breve B_{R}^{+}\,(=\{(r,y):r> 0,\,r^{2}+|y|^{2}<R^{2}\})$.
Direct computation then shows that
\[%
 \frac{1}{r^{1+\beta}}\frac{\partial}{\partial r}\Bigl(r^{1+\beta} \frac{\partial h}{\partial r}\Bigr) + %
\sum_{k=1}^{\ell}\frac{\partial^{2}h}{\partial y_{k}^{2}}=0  %
\dl{beta-lap} %
\]%
on $\breve B_{R}^{+}$, with
\[%
h(r,y)=h_{j}(r,y)=r^{-\gamma_{j}}\int_{\Sigma}v(r\omega,y)\,\varphi_{j}(\omega)\,d\mu  %
 \text{ and }\beta=\beta_{j}=2\sqrt{(\fr{n-2}{2})^{2}+\lambda_{j}}.
\dl{beta-lap-1} %
\]%
A solution $h$ of~\ref{beta-lap} (with any $\beta> 0$) which has the property
\[%
\int_{\{(r,y):r^{2}+|y|^{2}<R^{2}\}}r^{-2}h^{2}\,r^{1+\beta}\,drdy<\infty %
\dl{l2-harm-1}
\]%
will be referred to as a \emph{$\beta$-harmonic function on} $\{(r,y):r\ge 0,\, r^{2}+|y|^{2}<R^{2}\}$.  Using~\ref{l2-harm-1}
together with the weak form of the equation~\ref{beta-lap} this is easily shown to imply the $W^{1,2}$ estimate
\[%
\int_{\{(r,y):r^{2}+|y|^{2}<\rho^{2}\}}\bigl((D_{r}h)^{2}+  %
                           \tsum_{j=1}^{\ell}(D_{y_{j}}h)^{2}\bigr)\,r^{1+\beta}\,drdy<\infty %
\dl{l2-harm-2}
\]%
for each $\rho<R$, and there is a unique such function with prescribed $C^{1}$ data $\varphi$ on
$\bigl\{(r,y):r^{2}+|y|^{2}=R^{2}\bigr\}$, obtained, for example, by minimizing \smash{$\int_{r^{2}+|y|^{2}\le
R^{2}}(u_{r}^{2}+|u_{y}|^{2})\,r^{1+\beta}d\mu$} among functions with trace $\varphi$ on $r^{2}+|y|^{2}=R^{2}$.

If $\beta$ is an integer (e.g.\ $\beta=1$ when $n=7$, $j=1$, and $\C_{0}$ is the ``Simons cone''
$\bigl\{x\in\R^{8}:\sum_{i=1}^{4}(x^{i})^{2}=\sum_{i=5}^{8}(x^{i})^{2}\bigr\}$) then the $\beta\,$-Laplacian operator as in
\ref{beta-lap} is just the ordinary Laplacian in $\mathbb{R}^{2+\beta+\ell}$, at least as it applies to functions $u=u(r,y)$ with
$r=|x|$.  Even for $\beta$ non-integer, there is, analogous to the case when $\beta$ is an integer, $\beta$-harmonic polynomials of
each order (i.e.\ homogeneous polynomial solutions $h$ of \ref{beta-lap}) of each order $q=0,1,\ldots$.

Indeed, as shown in the Appendix below (extending the discussion of~\cite{Sim94} to arbitrary $\ell$ and at the same time
showing the relevant power series converge in the entire ball $\bigl\{ (r,y):r\ge 0,\,r^{2}+|y|^{2}<R^{2}\bigr\}$ rather than merely in
$\bigl\{ (r,y):r\ge 0,\,r^{2}+|y|^{2}<\theta R^{2}\bigr\}$ for suitable $\theta\in (0,1)$), if $h=h(r,y)$ is a solution of the weak
form of~\ref{beta-lap} on $\{(r,y):r> 0,\, r^{2}+|y|^{2}<R^{2}\}$, i.e.\ 
\[%
\int_{\breve B_{1}^{+}} (u_{r}\zeta_{r}+u_{y}\cdot \zeta_{y})\,d\mu_{+}=0
\dl{weak-form-0}
\]%
for all Lipschitz $\zeta$ with compact support in $r^{2}+|y|^{2}<R^{2}$ and if $u$ satisfies~\ref{l2-harm-1}, then $h$ is a real
analytic function of the variables $r^{2},y$, and, on all of $\{(r,y):r\ge 0,\, r^{2}+|y|^{2}<R^{2}\}$, $h$ is a convergent sum
\[%
h(r,y)=\tsum_{q=0}^{\infty}h_{q}(r,y) \,\,\,\forall (r,y) \text{ with }r\ge 0 \text{ and } \sqrt{r^{2}+|y|^{2}}<R, %
\dl{h-exp}
\]%
where $h_{q}$ is the degree $q$ homogeneous $\beta$-harmonic polynomial in the variables $r^{2},y$ obtained by
selecting the order $q$ terms in the power series expansion of $h$.

In case $\ell=1$ there is, up to scaling, a unique $\beta$-harmonic polynomial of degree $q$ of the form
{\belowdisplayskip8pt%
\[%
h_{q}=y^{q}+\sum_{k>0,\ell\ge 0,2k+j=q}c_{k\ell}r^{2k}y^{j}, %
\]%
By direct computation, in case $\ell=1$ the $\beta$-harmonic polynomials $h_{0},h_{1},h_{2},h_{3},h_{4}$ are
respectively 
\[%
1,\,\,\,y,\,\,\,y^{2}-\fr{1}{2+\beta}r^{2},\,\,\,y^{3}-\fr{3}{2+\beta}r^{2}y,\,\,\, %
y^{4}-\fr{6}{2+\beta}r^{2}y^{2}+\fr{3}{(2+\beta)(4+\beta)}r^{4}. %
\]}%
In the case $\ell\ge 2$, for each homogeneous degree $q$ polynomial $p_{0}(y)$ there is a unique $\beta$-harmonic
homogeneous polynomial $h$ of degree $q$, with
\[%
h(r,y)=
\begin{cases}%
p_{0}(y)+\tsum_{j=1}^{q/2}r^{2j}p_{j}(y), &q \text{ even } \\ %
p_{0}(y)+\tsum_{j=1}^{(q-1)/2} r^{2j}p_{j}(y),& q\text{ odd}, 
\end{cases}%
\]%
where each $p_{j}(y)$ is a homogeneous degree $q-2j$ polynomial, and the $p_{j}$ are defined inductively by
$p_{j+1}=-(2j+2)^{-1}(2j+2+\beta)^{-1}\Delta_{y}p_{j}$, $j\ge 0$.

In terms of spherical coordinates $\rho=\sqrt{r^{2}+|y|^{2}}$ and $\omega=\rho^{-1}(r,y) \in \Sph^{\ell}_{+}$,
\[%
\Sph^{\ell}_{+}=\{\omega=(\omega_{1},\ldots,\omega_{\ell})\in\R^{\ell+1}:\omega_{1}>0,\,|\omega|=1\}, 
\]%
the $\beta$-Laplacian $\Delta_{\beta}$ (i.e.\ the operator on the left of~\ref{beta-lap}) is
 \[%
    \Delta_{\beta}h = \rho^{-\ell-1-\beta}\frac{\partial}{\partial \rho} %
   \Bigl(\rho^{\ell+1+\beta}\frac{\partial h}{\partial \rho}\Bigr)+ \rho^{-2}\omega_{1}^{-1-\beta} %
\dvg_{\Sph_{+}^{\ell}}\bigl(\omega_{1}^{1+\beta}\nabla_{\Sph_{+}^{\ell}}h\bigr). %
\dl{sph-coords}
\]%
Notice that in the  case $\ell=1$ we can write $h=h(r,\theta)$ where $\theta=\arctan (y/r)\in
(-\pi/2,\pi/2)$ and~\ref{sph-coords}  can be written
 \[%
    \Delta_{\beta}h = \rho^{-2-\beta}\frac{\partial}{\partial \rho} %
   \Bigl(\rho^{2+\beta}\frac{\partial h}{\partial \rho}\Bigr)+ \rho^{-2}\cos^{-1-\beta}\theta
\frac{\partial}{\partial\theta}\bigl(\cos^{1+\beta}\theta\frac{\partial h}{\partial \theta}\bigr).
\]%
Using~\ref{sph-coords} we see that the order $q$ homogeneous $\beta$-harmonic polynomials $h_{q}$ (which are
homogeneous of degree $q$ in the variable $\rho$) satisfy the identity
{\abovedisplayskip8pt\belowdisplayskip8pt%
\[%
\dvg_{\Sph_{+}^{\ell}}\bigl(\omega_{1}^{1+\beta}\nabla_{\Sph_{+}^{\ell}}h_{q}\bigr)  %
= q(q+\ell+\beta)h_{q} \omega_{1}^{1+\beta}. %
\dl{jth-eig}
\]%
Hence we have the orthogonality of $h_{p},h_{q}$ for $p\neq q$ on
\[%
\Sph^{\ell}_{+}=\{(r,y):r> 0,\,r^{2}+|y|^{2}=1\} %
\dl{Sell} %
\]}%
with respect to the measure $d\nu_{+}=\smash{{\omega_{1}}^{\! 1+\beta}}d\mu$ ($\mu=\ell$-dimensional
Hausdorff measure on $\Sph_{+}^{\ell}$):
\[%
\int_{\Sph_{+}^{\ell}}h_{p}\,h_{q}\,d\nu_{+}=0 \text{ for $p\neq q$},\quad d\nu_{+}=\omega_{1}^{1+\beta}d\mu. %
\dl{orthog}
\]%
Thus if $h$ satisfies~\ref{beta-lap} and~\ref{l2-harm-1} then 
\begin{align*}%
  \int_{B^{+}_{R}}h^{2}(r,y)\,r^{1+\beta}drdy   %
  & =\int_{0}^{R}\int_{\Sph_{+}^{\ell}}h^{2}(\rho\omega)\, %
    \rho^{1+\beta}\,d\nu_{+} \,\rho^{\ell} d\rho \dtg{l2-norm}\\ %
  \noalign{\vskip-3pt}
  &=\sum_{q=0}^{\infty}(\ell+2+\beta+2q)^{-1}N_{q}^{2}R^{\ell+2+\beta+2q}, %
\end{align*}%
where $B_{R}^{+}=\{(r,y):r\ge 0,\,\,r^{2}+|y|^{2}\le R\}$ and
$N_{q}=\bigl(\int_{\Sph_{+}^{\ell}}h_{q}^{2}(\omega)\,\omega_{1}^{1+\beta} d\mu\bigr)^{1/2}$ and $h_{q}(r,y)$ are as
in~\ref{h-exp}.

Using~\ref{orthog} it is shown in the Appendix that the homogeneous $\beta$-harmonic polynomials are complete in
$L^{2}(\nu_{+})$, where $d\nu_{+}=\omega_{1}^{1+\beta}d\mu_{\ell}$ on $\Sph^{\ell}_{+}$.  Thus each $\varphi\in
L^{2}(\nu_{+})$ can be written as an $L^{2}(\nu_{+})$ convergent series
{\abovedisplayskip8pt\belowdisplayskip8pt%
\[%
\varphi=\sum_{q=0}^{\infty}h_{q}|\Sph^{\ell}_{+}, %
\dl{b-complete}
\]%
with each $h_{q}$ either zero or a homogeneous degree $q$ $\beta$-harmonic polynomial.

Next observe that if
\[%
\mathcal{L}_{\C}v=0 \text{ on }\C\cap B_{R} \text{ with } \int_{\C\cap B_{R}}|x|^{-2}v^{2}\,d\mu<\infty, %
\dl{proper-soln-LC}
\]}%
then each $h_{j}(r,y)$ defined as in~\ref{beta-lap-1} does satisfy the condition~\ref{l2-harm-1}, and hence we have
\[%
v=\sum_{j=1}^{\infty}\sum_{q=0}^{\infty}r^{\gamma_{j}}h_{j,q}(r,y)\varphi_{j}(\omega), %
\dl{LC-soln}
\]%
where $h_{j,q}$ is a homogeneous degree $q$ $\beta_{j}$-harmonic polynomial and, using the orthogonality~\ref{orthog},
\[%
\int_{\C\cap B_{\rho}}v^{2}\,d\mu= \rho^{\ell+2+\beta_{1}} \sum_{j=1}^{\infty}\sum_{q=0}^{\infty} %
(2+\ell+\beta_{j}+2q)^{-1}N_{j,q}^{2}\rho^{2q+\beta_{j}-\beta_{1}}, \,\,\rho<R,%
\dl{l2-norm-LC}
\]%
where $N_{j,q}=\bigl(\int_{\Sph_{+}^{\ell}}h_{j,q}^{2}(\omega)\, d\nu_{+}\bigr)^{1/2}$.

Observe that (except for $(j,q)=(1,0)$) it is possible that $r^{\gamma_{j}}h_{j,q}$ could have the same homogeneity as
$r^{\gamma_{i}}h_{i,p}$ (i.e.\ $\gamma_{j}+q=\gamma_{i}+p$) for some $i\neq j,\,p\neq q$, but in any case, after some
rearrangement of the terms in~\ref{l2-norm-LC}, we see that there are exponents
{\abovedisplayskip8pt\belowdisplayskip8pt%
\[%
0=q_{1}<q_{2}< q_{2}< \cdots < q_{i}\to\infty,\quad q_{i}=q_{i}(\C_{0}),
\]%
($q_{i}$ not necessarily integers, but nevertheless fixed depending only on $\C_{0}$) such that if $v$ is as
in~\ref{proper-soln-LC} then
\[%
\avint_{\C\cap B_{\rho}}v^{2} = \sum_{i=1}^{\infty} b_{i}^{2}\rho^{2q_{i}} 
\dl{l2-norm-v}
\]%
(which in particular is an increasing function of $\rho$), $\rho\in (0,R]$, for suitable constants $b_{j},\,j=1,2,\ldots$, where
we use the notation
\[%
\avint_{\C\cap B_{\rho}}f = \rho^{-\ell-2-\beta_{1}}\int_{\C\cap B_{\rho}}f\,d\mu, \quad \rho>0. %
\dl{def-avint}
\]%
We claim that the logarithm of the right side of~\ref{l2-norm-v} is a convex function $\psi(t)$ of $t=\log R$:
\[%
\psi''(t)\ge 0\text{ where }%
\psi(t)=\log\bigl(\avint_{\C\cap B_{\rho}}v^{2}\bigr|_{\rho=e^{t}}\bigr)\,\,\bigl(\,=\log(\tsum_{i}b_{i}^{2}e^{2q_{i}t}) %
\text{ by~\ref{l2-norm-v}}\bigr). %
\dl{conv-psi}
\]%
To check this, note that 
\[%
\psi''(t) =\psi^{-2}(t)\bigl((\tsum_{i}b_{i}^{2}e^{2q_{i}t})(\tsum_{i}4q_{i}^{2}b_{i}^{2}e^{2q_{i}t})- %
(\tsum_{i}2q_{i}b_{i}^{2}e^{2q_{i}t})^{2} \bigr)
\]%
and by Cauchy-Schwarz this is non-negative for $t\in (0,R)$ and if there is a $t_{0}\in (-\infty,\log R)$ such that
$\psi''(t_{0})=0$ then there is $i_{0}\in\{1,2,\ldots\}$ such that $b_{i}=0$ for every $i\neq i_{0}$, in which case
\[%
\psi(t)=\log b_{i_{0}}+2q_{i_{0}}t \quad t\in (-\infty,\log R). %
\dl{psi-lin}
\]%
In particular the convexity~\ref{conv-psi} of $\psi$ implies
\[%
\psi(t) -\psi(t-\log 2)\ge \psi(t-\log 2) -\psi(t-2\log 2),\quad t\in (-\infty,\log R),%
\dl{eq-rat}
\]%
and equality at any point $t\in (-\infty,\log R)$ implies~\ref{psi-lin} holds for some $i_{0}$ and the common value of each side
of~\ref{eq-rat} is $\log 4^{q_{i_{0}}}$.  Thus we see that if $Q\in (0,\infty)\setminus\{4^{q_{1}},4^{q_{2}},\ldots\}$ and $v$ is
not identically zero then for each given $\rho\in (0,R]$
\[%
\avint_{\C\cap B_{\rho/2}}v^{2}\bigg/\avint_{\C\cap B_{\rho/4}}v^{2}\ge Q\Longrightarrow %
\avint_{\C\cap B_{\rho}}v^{2}\bigg/\avint_{\C\cap B_{\rho/2}}v^{2}>Q.
\dl{lin-growth}
\]}%

\noindent{\bf{}\tl{rem-lin-growth} Remark: }
Notice that if $v_{1},\ldots,v_{q}$ are solutions of $\mathcal{L}_{\C}v=0$ on $\C\cap B_{R}$ satisfying
\[%
\avint_{\C\cap B_{R}}|x|^{-2}v_{j}^{2}(x,y)\,d\mu(x,y)<\infty \text{ for each }j=1,\ldots,q 
\]%
with $v_{j}\neq 0$ for some $j$, then $\smavint_{\C\cap B_{R}}\smash{\sum_{j=1}^{q}v_{j}^{2}\,d\mu}$ has again the
form of~\emph{\ref{l2-norm-v}}, so the implication~\emph{\ref{lin-growth}} applies with the sum $\sum_{j=1}^{q}v_{j}^{\,2}$
in place of~$v^{2}$.

\section[Growth estimates for solutions of $\smash{\mathcal{L}_{\!M}w=0}$]{Growth estimates for solutions
of $\mathcal{L}_{\!M}w=0$}\label{growth}

\noindent Here we discuss growth estimates analogous to those used in~\cite{Sim83a}, \cite{Sim85} for solutions of
$\mathcal{L}_{M}w=0$; in particular we discuss conditions on $M$ and $w$ which ensure that the growth behaviour of
solutions $w$ of $\mathcal{L}_{M}w=0$ is analogous to that of the solutions $v$ of $\mathcal{L}_{\C}v=0$ discussed in the
previous section.

The main growth lemma below (Lemma~\ref{growth-lem}) applies locally in balls (so $M$ could be a complete minimal
hypersurface in a ball $\breve B_{R}$ rather than the whole space), and we do not need the inclusion $M\subset
U_{+}\times\R^{\ell}$. We in fact assume $R,\lambda,\Lambda>0$ and
{\abovedisplayskip7pt\belowdisplayskip7pt%
\[%
\left\{\begin{aligned}%
&\text{\,$M\subset \breve B_{R}$ is embedded, minimal, and satisfies~\ref{str-stab-M}}  %
\,\,\text{ for every }\,\zeta\in C^{1}_{c}(\breve B_{R}), \\  %
\noalign{\vskip-3pt} %
& \, R^{-n-\ell-2}\int_{M\cap B_{R}}d^{2}(x)\,d\mu(x,y)<\epsilon, \text{ and } %
                                                                             R^{-n-\ell}\mu(M\cap B_{R})\le \Lambda, %
 \end{aligned}\right.%
\dl{local-M}
 \]%
 with an $\epsilon$ (small) to be specified and with 
\[%
d(x)=\dist(x,\C_{0})\,\,(=\dist((x,y),\C))\text{ for }(x,y)\in\R^{n+1}\times\R^{\ell}. %
\dl{def-d}
\]%
Taking $\zeta$ in~\ref{w-1-2-est} to be $1$ in $B_{\theta R}$, $0$ on $\partial B_{R}$, and $|D\zeta|\le 2/(1-\theta)$, we
obtain
\[%
\int_{M\cap B_{\theta R}}\hskip-5pt |x|^{-2}w^{2}(x,y)\,d\mu(x,y)\le  %
CR^{-2}\int_{M\cap B_{R}}\hskip-5pt w^{2}\,d\mu\,\,\forall\theta\in [\ha,1), \,\, C=C(\lambda,\theta). %
\dl{pre-gro-1}
\]%
Notice that  if we have a ``doubling condition''
\[%
\int_{M\cap B_{R}}w^{2}(x,y)\,d\mu(x,y)\le  K\int_{M\cap B_{R/2}}w^{2}\,d\mu,
\dl{doubling}
\]%
then \ref{pre-gro-1} implies that
\[%
\int_{M\cap B_{\theta R}}\hskip-5pt |x|^{-2}w^{2}(x,y)\,d\mu(x,y)\le CKR^{-2}\int_{M\cap B_{R/2}}w^{2}\,d\mu,
\]%
so, for $\delta^{-2}\ge 2CK$ (i.e.\ $\delta\le (2CK)^{-1/2}$) we have
\[%
\int_{M\cap B_{\theta R}}\hskip-5pt |x|^{-2}w^{2}(x,y)\,d\mu(x,y)\le  %
CR^{-2}\int_{M\cap B_{R/2}\setminus\{(x,y):|x|\le \delta R\}}\hskip-5pt w^{2}\,d\mu\,\,\forall\theta\in [\ha,1),%
\dl{pre-gro-2}
\]}%
where $C=C(\lambda,\theta,K)$.  The following lemma shows that in fact the inequality~\ref{pre-gro-2} can be improved.

\begin{state} {\bf{}\tl{improved-doubling} Lemma.}%
For each $\lambda,\Lambda,K>0$ and $\theta\in \smash{[\ha,1)}$ there is $\epsilon=\epsilon(\lambda,\Lambda,\theta,
K,\C)>0$ such that if $M$ satisfies~\emph{\ref{local-M}} and $\mathcal{L}_{M}w=0$ on $M\cap\breve B_{R}$, and if the
doubling condition~\emph{\ref{doubling}} holds then
\[%
\int_{M\cap B_{\theta R}}|x|^{-2}w^{2}(x,y)\,d\mu(x,y)\le  %
CR^{-2}\int_{M\cap B_{R/2}\setminus\{(x,y):|x|\le \frac{1}{3} R\}}w^{2}\,d\mu,  %
\]%
where $C=C(\lambda,\Lambda,\theta,K,\C)$.
\end{state}%

\begin{proof}{\bf{}Proof:} Let $\theta\in [\fr{3}{4},1)$.  By scaling we can assume $R=1$.  If there is no $\epsilon$
ensuring the claim of the lemma, then there is a sequence $M_{k}$ of strictly stable hypersurfaces (with fixed $\lambda$)
with $\int_{M_{k}\cap B_{1}}d^{2}(x)\,d\mu(x,y)\to 0$ and a sequence $w_{k}$ of solutions of
$\mathcal{L}_{M_{k}}w_{k}=0$ such that
\[%
\int_{M_{k}\cap B_{1}}w_{k}^{2}(x,y)\,d\mu(x,y)\le  K\int_{M_{k}\cap B_{1/2}}w_{k}^{2}\,d\mu,
\pdl{doubling-k}
\]%
yet such that
\[%
\int_{M_{k}\cap B_{\theta}}|x|^{-2}w_{k}^{2}(x,y)\,d\mu(x,y)\ge  %
k\, \int_{M_{k}\cap B_{1/2}\setminus\{(x,y):|x|\le \frac{1}{3} \}}w_{k}^{2}\,d\mu.
\pdl{cont}
\]%
By~\ref{pre-gro-2}, 
\[%
\int_{M_{k}\cap B_{\theta}}|x|^{-2}w_{k}^{2}(x,y)\,d\mu(x,y)\le  %
C \int_{M_{k}\cap B_{1/2}\setminus\{(x,y):|x|\le \delta \}}w_{k}^{2}\,d\mu,
\pdl{pre-gro-4}
\]%
with $\delta=\delta(\lambda,\theta,K)$ and $C=C(\lambda,\theta,K)$. 

Since $\int_{M_{k}\cap B_{1}}d^{2}\,d\mu_{k}\to 0$, in view of~\ref{tangent-cone}\,(i),\,(ii) there are sequences
$\eta_{k}\downarrow 0$ and $\tau_{k}\uparrow 1$ with
\[%
\left\{\begin{aligned}%
&\,\,M_{k}\cap \breve B_{\tau_{k}}\setminus \{(x,y):|x|\le \eta_{k}\}\subset\cup_{j=1}^{q}\graph u_{k,j} \subset M_{k}\\ %
 &\qquad \sup(|u_{k,j}|+|\nabla_{\C}u_{k,j}|+|\nabla^{2}_{\C}u_{k,j}|) \to 0,  %
\end{aligned}\right.%
\pdl{gro-2a}
\]%
where $q=q(\Lambda)\ge 1$ and each $u_{k,j}$ is $C^{2}$ on a domain containing $\C\cap \breve B_{\tau_{k}}\setminus
\{(x,y):|x|\le \eta_{k}\}$.  Thus, with $\nu_{\C_{0}}$ the unit normal of $\C_{0}$ pointing into $U_{+}$ and with
\[%
 w_{k,j}(x,y)=w_{k}((x+u_{k,j}(x,y)\nu_{\C_{0}}(x),y)),  %
\pdl{gro-2c}
\]%
for $(x,y)\in \C\cap \breve B_{\tau_{k}} \setminus \{(x,y):|x|\le \eta_{k}\}$ and $j=1,\ldots,q$, we see that $w_{k,j}$ satisfies a
uniformly elliptic equation with coefficients converging in $C^{1}$ to the coefficients of $\mathcal{L}_{\C}$ on
$\Omega_{\sigma,\delta}=\C\cap \breve B_{\sigma}\setminus\{(x,y):|x|>\delta\}$ for each $\sigma\in [\ha,1)$ and $\delta\in
(0,\ha]$, hence, by Schauder estimates and~\ref{doubling-k}, 
\[%
|w_{k,j}|_{C^{2}(\Omega_{\sigma,\delta})}\le C\|w_{k}\|_{L^{2}(B_{1})}\le CK\|w_{k}\|_{L^{2}(B_{1/2})} 
\pdl{c2-ests}
\]%
with $C$ independent of $k$, $k\ge k(\sigma,\delta)$. Hence, by~\ref{gro-2a} and~\ref{c2-ests},
\[%
\tilde w_{k,j}=(\tint_{M_{k}\cap B_{1/2}} w_{k}^{2}\,d\mu_{k})^{-1/2} w_{k,j}  %
\]%
has a subsequence converging in $C^{1}$ locally in $\C\cap\breve B_{1}$ to a smooth solution $v_{j}$ of
$\mathcal{L}_{\C}v_{j}=0$ on $\C\cap\breve B_{1}$, and since by~\ref{pre-gro-4}  $\int_{M_{k}\cap B_{1/2}\cap
\{(x,y):|x|\le \sigma\}}w_{k}^{2}\le C\sigma^{2}\int_{M_{k}\cap B_{1/2}}w_{k}^{2}$ for all $\sigma\in (0,\ha]$, 
we can then conclude
{\abovedisplayskip8pt\belowdisplayskip8pt%
\[%
\int_{\C\cap B_{1/2}}\sum_{j=1}^{q}v_{j}^{2}\,d\mu=1. %
\pdl{norm-1}
\]%
But by~\ref{doubling-k} and \ref{cont} 
\[%
\int_{M_{k}\cap B_{1/2}\setminus\{(x,y):|x|\le \frac{1}{3} \}}w_{k}^{2}\,d\mu\le 
Ck^{-1} \int_{M_{k}\cap B_{1/2}}w_{k}^{2}\,d\mu,
\]%
and so, multiplying each side by $(\tint_{M_{k}\cap B_{1/2}} w_{k}^{2}\,d\mu_{k})^{-1/2}$ and taking limits, we
conclude
\[%
\int_{\C\cap B_{1/2}\setminus\{(x,y):|x|\le \frac{1}{3}\}} \sum_{j=1}^{p} v_{j}^{2} = 0.
\]}%
In view of~\ref{norm-1} this contradicts unique continuation for solutions of $\mathcal{L}_{\C}v=0$ (applicable since the
solutions $v$ of $\mathcal{L}_{\C}v=0$ are real-analytic).\nobreak\end{proof}

\smallskip

We can now establish the growth lemma.

\begin{state}{\bf \tl{growth-lem} Lemma.}%
For each $\lambda,\Lambda>0$, $Q\in (0,\infty)\setminus\{4^{q_{1}},4^{q_{2}},\ldots\}$ ($q_{1},q_{2},\ldots$ as
in~{\rm\ref{l2-norm-v}}), and $\alpha\in [\ha,1)$, there is $\epsilon=\epsilon(Q,\alpha,\lambda,\Lambda,\C)>0$ such that if
$M$ satisfies~\emph{\ref{local-M}} and if $\mathcal{L}_{M}w=0$ on $M\cap\breve B_{R}$ then
\[%
\avint_{M\cap B_{R/2}}w^{2}\ge Q\avint_{M\cap B_{R/4}}w^{2}\Longrightarrow %
\avint_{B_{R}}w^{2} \ge Q\avint_{B_{R/2}}w^{2},  %
\leqno{\rm (i)} %
\]%
and
\[%
 \avint_{M\cap B_{R}}w^{2} \ge \alpha \avint_{M\cap B_{R/2}}w^{2},  %
\leqno{\rm (ii)}
\]%
where we use the notation (analogous to~{\rm\ref{def-avint}})
\[%
\avint_{M\cap B_{\rho}}f= \rho^{-\ell-2-\beta}\int_{M\cap B_{\rho}}f\,d\mu.
\]%
Notice that no hypothesis like \smash{$M\subset U_{+}\times\R^{\ell}$} is needed here.
\end{state}%

\begin{proof}{\bf{}Proof:} By scaling we can assume $R=1$.  If there is $Q\in
(0,\infty)\setminus\{4^{q_{1}},4^{q_{2}},\ldots\}$ such that there is no $\epsilon$ ensuring the first claim of the lemma,
then there is a sequence $M_{k}$ satisfying \ref{str-stab-M} (with fixed $\lambda$) with $\int_{M_{k}\cap
B_{1}}d^{2}(x)\,d\mu(x,y)\to 0$ and a sequence $w_{k}$ of
solutions of $\mathcal{L}_{M_{k}}w_{k}=0$, such that
{\abovedisplayskip8pt\belowdisplayskip8pt%
\[%
\avint_{M_{k}\cap B_{1/2}}w_{k}^{2}\ge Q\avint_{M_{k}\cap B_{1/4}}w_{k}^{2}\text{ and }
\avint_{M_{k}\cap B_{1}}w_{k}^{2}<Q\avint_{M_{k}\cap B_{1/2}}w_{k}^{2}. %
\pdl{gro-2}
\]}%
The latter inequality implies we have the doubling condition~\ref{doubling} with $K=Q$, so we can repeat the compactness
argument in the proof Lemma~\ref{improved-doubling}. Thus by~\ref{gro-2}, with the same notation as in the proof
of~\ref{improved-doubling}, we get convergence of $w_{k,j}$ to a smooth solution $v_{j}$ of $\mathcal{L}_{\C}v_{j}=0$ on
$\C\cap\breve B_{1}$ with $\smavint_{\C\cap B_{1/2}}\tsum_{j=1}^{q}v_{j}^{2}\,d\mu=1$, $\int_{\C\cap
B_{1}}|x|^{-2}\tsum_{j=1}^{q}v_{j}^{2}<\infty$,
\begin{align*}%
&\avint_{\C\cap B_{1/2}}\tsum_{j=1}^{q}v_{j}^{2}\,d\mu\ge Q\avint_{\C\cap B_{1/4}}\tsum_{j=1}^{q}v_{j}^{2}\,d\mu %
\text{  and }\ptg{gro-3} \\ %
&\hskip1in\avint_{\C\cap B_{1}}\tsum_{j=1}^{q}v_{j}^{2}\,d\mu\le Q\avint_{\C\cap B_{1/2}}\tsum_{j=1}^{q}v_{j}^{2}\,d\mu.
 \end{align*}%
In view of~\ref{rem-lin-growth}, this contradicts~\ref{lin-growth}.

Similarly, if the second claim of the lemma fails for some $\alpha\in [\ha,1)$, after taking a subsequence of $k$ (still
denoted $k$) we get sequences $M_{k}, w_{k}$ and $w_{k,j}$  with
\[%
\avint_{M_{k}\cap B_{1}} w_{k}^{2}\,d\mu<\alpha\avint_{M_{k}\cap B_{1/2}} w_{k}^{2}\,d\mu
\]%
(i.e.\ the doubling condition~\ref{doubling} with $K=\alpha$), and smooth solutions $v_{j}=\lim w_{k,j}$ of
$\mathcal{L}_{\C}v_{j}=0$ on $\C\cap\breve B_{1}$ with $0<\int_{M\cap B_{1}}|x|^{-2}v_{j}^{2}<\infty$ and
\[%
\avint_{\C\cap B_{1}}\tsum_{j=1}^{q}v_{j}^{2}\,d\mu\le\alpha\avint_{\C\cap B_{1/2}}\tsum_{j=1}^{q}v_{j}^{2}\,d\mu,
\]%
which is impossible by~\ref{l2-norm-v} and~\ref{rem-lin-growth}.
\end{proof}

\smallskip

Since $\mathcal{L}_{M}$ is the linearization of the minimal surface operator $\mathcal{M}_{M}$, a smooth family
$\{M_{t}\}_{|t|<\epsilon}$ of minimal submanifolds with $M_{0}=M$ must generate a velocity vector $v$ at $t=0$ with
normal part $w=v\cdot\nu_{M}$ ($\nu_{M}$ a smooth unit normal of $M$) being a solution of the equation
$\mathcal{L}_{M}w=0$.  In particular the family of homotheties $(1+t)M$ generates the solution $w=(x,y)\cdot
\nu_{M}(x,y)$ and the translates $M+t e_{n+1+j}$ generate the solutions $w=e_{n+1+j}\cdot\nu_{M}(x,y) =
\nu_{y_{j}}(x,y)$, $j=1,\ldots,\ell$. Thus
\[%
w=(x,y)\cdot\nu_{M} \text{ and } w=\nu_{y_{j}} \text{ both satisfy }\mathcal{L}_{M}w=0\text{ and the %
  inequality }\ref{pre-gro-1}.  %
\dl{jac-solns}
\]%

\smallskip

\begin{state}{\bf \tl{strong-doub-lem} Corollary.}%
Suppose $\Lambda,\lambda,\gamma>0$, $q\in \{1,2,\ldots\}$, and assume $M$ is a complete, embedded, minimal
hypersurface in $\R^{n+1+\ell}$, strict stability~{\rm\ref{str-stab-M}} holds, and $M$ has $\C$ with some multiplicity $q$ as
its unique tangent cone at $\infty$.  If $\mathcal{L}_{M}w=0$, and $\sup_{R>1}R^{-\gamma}\smavint_{M\cap
B_{R}}w^{2}<\Lambda$, then there is $R_{0}=R_{0}(M,\gamma,\lambda,\Lambda,q)$ such that for all $R\ge R_{0}$ we
have the ``strong doubling condition''
\[%
\avint_{M\cap B_{R}}|x|^{-2}w^{2}\le CR^{-2}\avint_{M\cap B_{R/2}\setminus\{(x,y):|x|<\smfr{1}{3} R\}}w^{2}, %
\quad C =C(\gamma,q,\lambda,\Lambda). %
\]%
\end{state}%

\noindent{\bf \tl{strong-doub-rem} Remarks:} {\bf (1)} Observe that in case $M\subset U_{+}\times\R^{\ell}$ the unique
tangent cone assumption here is automatically satisfied by virtue of {\rm Lemma~\ref{tangent-cone}}.

{\bf{}(2)} Since $\mu(M\cap B_{R})\le q\mu(\C\cap B_{1})R^{n+\ell}$, $|(x,y)\cdot\nu_{M}(x,y)|\le R$, and
$|\nu_{y_{j}}|\le 1$ on $M\cap B_{R}$, in view of~{\rm\ref{jac-solns}} the above lemma applies to both
$w=(x,y)\cdot\nu_{M}$ and $w=\nu_{y_{j}}$ with $\gamma=2|\gamma_{1}|+2$ and $\gamma=2|\gamma_{1}|$
respectively.

\medskip

\begin{proof}{\bf{}Proof of Corollary~\ref{strong-doub-lem}:} We can of course assume that $w$ is not identically zero.  In view of
Lemma~\ref{tangent-cone}, for any $\tilde\gamma>\gamma$ with $2^{\tilde\gamma}\in
(1,\infty)\setminus\{4^{q_{1}},4^{q_{2}},\ldots\}$, $M\cap\breve B_{R}$ satisfies the hypotheses of Lemma~\ref{growth-lem} with
$Q=2^{\tilde\gamma}$ for all sufficiently large $R$, so
{\abovedisplayskip8pt\belowdisplayskip8pt%
\[%
\avint_{M\cap B_{2^{k}R}}w^{2}\ge 2^{\tilde\gamma} \avint_{M\cap B_{2^{k-1}R}}w^{2}\Longrightarrow %
\avint_{M\cap B_{2^{k+1}R}}w^{2}\ge 2^{\tilde \gamma}\avint_{M\cap B_{2^{k}R}}w^{2} %
\]%
for $k=1,2,\ldots$ and for any choice of $R$ sufficiently large (depending on $M,\tilde\gamma,\lambda$), and hence, by
iteration, if $\smavint_{M\cap B_{R}}w^{2}\ge 2^{\tilde\gamma} \smavint_{M\cap B_{R/2}}w^{2}$ we would have
$\smavint_{M\cap B_{2^{k}R}}w^{2}\ge C2^{k\tilde\gamma}$, $k=1,2,\ldots$ with $C>0$ independent of $k$ ,
contrary to the hypothesis on $w$ since $C2^{k(\tilde\gamma-\gamma)}>\Lambda$ for sufficiently large $k$. 
Thus, with such a choice of $\tilde\gamma$, we have the doubling condition~\ref{doubling} with $K=2^{\tilde\gamma}$:
\[%
\avint_{M\cap B_{R}}w^{2}\le 2^{\tilde\gamma}\avint_{M\cap B_{R/2}}w^{2}
\]}%
for all $R\ge R_{0}$, $R_{0}=R_{0}(M,\tilde\gamma)$.  Then the required result is proved by
Lemma~\ref{improved-doubling}.\nobreak\end{proof}

\section[Asymptotic $\smash{L^{2}}$ Estimates for $d$]{Asymptotic $L^{2}$ Estimates
for $d$} \label{dist-fn}

\noindent The following lemma gives an $L^{2}$ estimate for the distance function $d|M\cap B_{R}$ as $R\to \infty$ in case
$M$ is as in Theorem~\ref{main-th}, where as usual
{\abovedisplayskip8pt\belowdisplayskip8pt%
\[%
d(x)=\dist(x,\C_{0})\,\,(=\dist((x,y),\C)\text{ for } (x,y)\in\R^{n+1}\times\R^{\ell}).
\]%
In this section we continue to use the notation
\[%
\avint_{M\cap B_{\rho}}f =\rho^{-\ell-2-\beta_{1}}\int_{M\cap B_{\rho}}f\,d\mu.
\]}%

\begin{state}{\bf\tl{main-dist-est} Lemma.}%
Let $\alpha,\lambda,\Lambda>0$ and assume $M$ is embedded, minimal, $M\subset U_{+}\times\R^{\ell}$,
$R^{-n-\ell}\mu(M\cap B_{R})\le \Lambda$ for each $R>1$ and $M$ satisfies the strict stability~{\rm\ref{str-stab-M}}. Then there
is $R_{0}=R_{0}(M,\lambda,\Lambda,\alpha)>1$ such that
\[%
C^{-1}R^{-\alpha}\avint_{M\cap B_{R_{0}}}d^{2}\le \avint_{M\cap B_{R}}d^{2} \le CR^{\alpha}\avint_{M\cap
B_{R_{0}}}d^{2} \,\,\,\,\forall R \ge R_{0}, \quad C=C(\lambda,\Lambda,\C_{0}).
\]%
\end{state}%

To facilitate the proof we need the technical preliminaries of the following lemma.  In this lemma 
\[%
\nu_{M}=(\nu_{1},\ldots,\nu_{n+\ell+1})
\dl{normal-M}
\]%
continues to denote a smooth unit normal for $M$, $\nu_{\C_{0}}$ continues to denote the unit normal of $\C_{0}$
pointing into $U_{+}$, and $\epsilon_{0}=\epsilon_{0}(\C_{0})\in (0,\ha]$ is assumed small enough to ensure that there is a
smooth nearest point projection
\[%
\pi:\bigl\{x\in\R^{n+1}:d(x)<\epsilon_{0}|x|\bigr\} \to \C_{0}. %
\dl{nearest-pt}
\]%

\begin{state}{\bf{}\tl{tech-prelim} Lemma.}%
Suppose $\Lambda>0$, $\delta\in (0,\ha]$, $\theta\in [\ha ,1)$,  $M$ satisfies~{\rm\ref{local-M}}, and $M\subset
U_{+}\times\R^{\ell}$.  Then there is $\epsilon=\epsilon(\delta,\lambda,\Lambda,\theta,\C)\in (0,\epsilon_{0}]$ such that, for
all $(x,y)\in M\cap B_{\theta R}$ with $d(x)<\epsilon|x|$,
\[%
 \min\bigl\{\bigl |\nu_{M}(x,y)-\bigl(\nu_{\C_{0}} (\pi(x)),0\bigr)\bigr|,  %
                                      \,\bigl|\nu_{M}(x,y)+\bigl(\nu_{\C_{0}}(\pi(x)),0\bigr)\bigr|\bigr\}<\delta  %
\leqno{\rm(i)}\\ %
\]%
in particular $|\nu_{y}|(x,y)\,\,(\,=|(\nu_{y_{1}}(x,y),\ldots,\nu_{y_{\ell}}(x,y))|)<\delta$), and 
\begin{align*}%
&(1-7\delta)|x|^{-1}|(x,0)\cdot \nu_{M}(x,y)|\le  \tg{\rm(ii)} \\ %
\noalign{\vskip-1pt}
&\hskip1in \bigl|(x,0)\cdot\nabla_{M}\bigl(d(x)/|x|\bigr)\bigr| \le  %
(1+7\delta)|x|^{-1}|(x,0)\cdot\nu_{M}(x,y)|.   %
\end{align*}%
\end{state}%

\noindent{\bf Remark:} The inequalities in~{\rm (ii)} do not depend on the minimality of $M$---{\rm (ii)} is true for any
smooth hypersurface $M$ at points where $d(x)<\epsilon|x|$ and where~{\rm (i)} holds for suitably small $\delta \in (0,\ha]$.

\begin{proof}{\bf{}Proof of Lemma~\ref{tech-prelim}:} By scaling we can assume $R=1$. If~(i) fails for some $\Lambda,
\delta,\theta$ then there is a sequence $(\xi_{k},\eta_{k})\in M\cap B_{\theta}$ with $d(\xi_{k})\le k^{-1}|\xi_{k}|$ and
\[%
\min\bigl\{|\nu(\xi_{k},\eta_{k})-\nu_{\C_{0}}(\pi(\xi_{k}))|,|\nu(\xi_{k},\eta_{k}) %
+\nu_{\C_{0}}(\pi(\xi_{k}))|\bigr\}\ge\delta. %
\pdl{close-0}
\]%
Then passing to a subsequence we have $\tilde\xi_{k}=|\xi_{k}|^{-1}\xi_{k}\to\xi\in\C_{0}$ and
$M_{k}=|\xi_{k}|^{-1}(M-(0,\eta_{k}))$ converges in the varifold sense to an integer multiplicity $V$ which is stationary in
$B_{\theta^{-1}}(\xi,0)$,  $\spt V\subset\overline U_{+}\times\R^{\ell}$, and $\xi\in\spt V\cap \C_{0}$. Also, by the
compactness and regularity theory of~\cite{SchS81}, $\sing V\cap \breve B_{\theta^{-1}}$ has Hausdorff dimension $\le
n+\ell-7$. So, since $\C_{0}$ is connected, by the maximum principle of~\cite{Ilm96} we have
\[%
V\res \breve B_{\theta^{-1}}(\xi,0)=(q\C+W)\res \breve B_{\theta^{-1}}(\xi,0) %
\pdl{close-1}
\]%
where $W$ is integer multiplicity, stationary in $\breve B_{\theta^{-1}}$ and with $\spt W\subset\overline
U_{+}\times\R^{\ell}$. Taking the maximum integer $q$ such that this holds, we then have
\[%
B_{\sigma}(\xi,0)\cap\spt W=\emptyset %
\pdl{close-2}
\]%
for some $\sigma>0$, because otherwise $\spt W\cap \overline\C\neq \emptyset$ and we could apply~\cite{Ilm96} again
to conclude $W\res\breve B_{\theta^{-1}}(\xi,0)=(\C+W_{1})\res\breve B_{\theta^{-1}}(\xi,0)$ where $W_{1}$ is
stationary integer multiplicity, contradicting the maximality of $q$ in~\ref{close-1}.

In view of~\ref{close-1}, \ref{close-2} we can then apply the sheeting theorem~\cite[Theorem 1]{SchS81} to assert that
$M_{k}\cap B_{\sigma/2}(\xi,0)$ is $C^{2}$ close to $\C$, and in particular
\[%
\min\bigl\{|\nu(\xi_{k},\eta_{k})-\nu_{\C_{0}}(\pi(\xi_{k}))|,|\nu(\xi_{k},\eta_{k}) +  %
\nu_{\C_{0}}(\pi(\xi_{k}))|\bigr\}\to 0,  %
\]%
contradicting~\ref{close-0}.

To prove~(ii), let $(x_{0},y_{0})\in M$ such that $d(x_{0})<\epsilon|x_{0}|$ and such that~(i) holds with
$(x,y)=(x_{0},y_{0})$, and let $\sigma>0$ be small enough to ensure that both $d(x)<\epsilon|x|$ and~(i) hold for all
$(x,y)\in M\cap B_{\sigma}(x_{0},y_{0})$.  Let
\[%
M_{y_{0}} = \bigl\{x\in\R^{n+1}:(x,y_{0})\in M\cap B_{\sigma}(x_{0},y_{0})\bigr\}.
\]%
Then, taking a smaller $\sigma$ if necessary, we can assume
\[%
M_{y_{0}} = \graph u=\bigl\{(\xi+u(\xi)\nu_{\C_{0}}(\xi),y_{0}):\xi\in\Omega\bigr\}, %
\pdl{close-4}
\]%
where $\nu_{\C_{0}}$ is the unit normal of $\C_{0}$ pointing into $U_{+}$, $\Omega$ is an open neighborhood of
$\pi(x_{0})$ in $\C_{0}$, and $u$ is $C^{2}(\Omega)$.

Then with $x=\xi+u(\xi)\nu_{\C_{0}}(\xi)\in M_{y_{0}}$ we have also $x(t)=t\xi+u(t\xi)\nu_{\C_{0}}(\xi)\in M_{ y_{0}}$
for $t$ sufficiently close to $1$ (because $\nu_{\C_{0}}(t\xi)=\nu_{\C_{0}}(\xi)$) and hence
\[%
\fr{d}{dt}\bigr|_{t=1}x(t)=\xi+\xi\cdot \nabla_{\C_{0}}u(\xi)\nu_{\C_{0}}(\xi)\in T_{x}M_{y_{0}}. %
\pdl{close-5}
\]%
With $\nu'(x)= (\nu_{1}(x,y_{0}),\ldots,\nu_{n+1}(x,y_{0}))$ (which is normal to $M_{y_{0}}$), we can assume, after
changing the sign of $\nu$ if necessary, that
\[%
\nu_{\C_{0}}(\xi)\cdot\nu'(x)>0 \text{ for }x\in M_{y_{0}}, \text{ so } |\nu_{\C_{0}}(\xi)-\nu'(x)|<\delta\text{ by (i)}.%
\pdl{close-5a}
\]%
By~\ref{close-5} we have
\[%
 x\cdot\nu'(x)  %
=\bigl(x-\fr{d}{dt}\bigr|_{t=1}x(t)\bigr)\cdot\nu'(x)=\bigl(u(\xi)-\xi\cdot\nabla_{\C_{0}}u(\xi)\bigr)\nu_{\C_{0}}(\xi)\cdot\nu'(x), %
\]%
hence
\[%
u(\xi)-\xi\cdot\nabla_{\C_{0}}u(\xi)= (\nu_{\C_{0}}(\xi)\cdot\nu'(x))^{-1}x\cdot\nu'(x). %
\pdl{close-6}
\]%
Differentiating the identity  $u(t\xi)=d\bigl(t\xi+u(t\xi)\nu_{\C_{0}}(\xi)\bigr)$ at $t=1$, we obtain
{\abovedisplayskip8pt\belowdisplayskip8pt%
\begin{align*}%
-\xi\cdot\nabla_{\C_{0}}u(\xi)  %
&= -\bigl(\xi+\xi\cdot\nabla_{\C_{0}}u(\xi)\nu_{\C_{0}}(\xi)\bigr)\cdot\nabla_{M_{y_{0}}}d(x)  \\ %
&=-\bigl(x-(u(\xi)-\xi\cdot\nabla_{\C_{0}}u(\xi))\nu_{\C_{0}}(\xi)\bigr)\cdot\nabla_{M_{y_{0}}}d(x) \\  %
 &=-x\cdot\nabla_{M_{y_{0}}}d(x) %
+\bigl((u(\xi)-\xi\cdot\nabla_{\C_{0}}u(\xi))\nu_{\C_{0}}(\xi)\bigr)\cdot\nabla_{M_{y_{0}}}d(x),   %
\end{align*}%
hence, adding $d(x)\,\,(\,=u(\xi))$ to each side of this identity, 
\begin{align*}%
d(x)-x\cdot\nabla_{M_{y_{0}}}d(x) %
&=\bigl(u-\xi\cdot \nabla_{\C_{0}}u(\xi)\bigr)\bigl(1-\nu_{\C_{0}}(\xi).\nabla_{M_{y_{0}}}d(x)\bigr)\ptg{close-7} \\ %
&=\bigl(\nu_{\C_{0}}(\xi)\cdot\nu'(x)\bigr)^{-1}x\cdot\nu'(x)
\bigl(1-\nu_{\C_{0}}(\xi).\nabla_{M_{y_{0}}}d(x)\bigr) %
\end{align*}}%
by~\ref{close-6}. Using the identity $|x|^{-1}x\cdot\nabla_{M_{y_{0}}}|x|=1-|\nu'(x)|^{-2}(x\cdot\nu'(x)/|x|)^{2}$ we see that
the left side of~\ref{close-7} is $\bigl(-x\cdot \nabla_{M_{y_{0}}}(d(x)/|x|)-|\nu'(x)|^{-2}(|x\cdot\nu'(x)|/|x|)^{2}\bigr)|x|$,
so~\ref{close-7} gives
{\abovedisplayskip8pt\belowdisplayskip0pt%
\[%
x\cdot \nabla_{M_{y_{0}}}(d(x)/|x|)=|x|^{-1} x\cdot\nu'(x)(1+E),  %
\]}%
where
{\abovedisplayskip6pt\belowdisplayskip8pt%
\[%
E=|\nu'(x)|^{-2}x\cdot\nu'(x)/|x|+(|\nu_{\C_{0}}(\xi)\cdot\nu'(x)|^{-1}-1)  %
-|\nu_{\C_{0}}(\xi)\cdot\nu'(x)|^{-1}\nu_{\C_{0}}(\xi)\cdot\nabla_{M_{y_{0}}}d(x). %
\]}%
Since $\nu_{\C_{0}}(\xi)\cdot\nabla_{M_{y_{0}}}d(x)=(\nu_{\C_{0}}(\xi)-\nu'(x))\cdot\nabla_{M_{y_{0}}}d(x)$ and
$x\cdot\nu'(x)= x\cdot\nu_{\C_{0}}(\xi) +x\cdot (\nu'(x)-\nu_{\C_{0}}(\xi))=
d(x)-x\cdot(\nu_{\C_{0}}(\xi)-\nu'(x))$, and $|\nu_{\C_{0}}(\xi)-\nu'(x)|<\delta$ by~\ref{close-5a}, direct
calculation then shows that $|E| < 7\delta$. Since $x\cdot \nabla_{M_{y_{0}}}(d(x)/|x|)=(x,0)\cdot
\nabla_{M}(d(x)/|x|)\bigr|_{(x,y_{0})}$ and $x\cdot\nu'(x)=(x,0)\cdot\nu(x,y_{0})$ we thus have the required
inequalities~(ii) at $(x_{0},y_{0})$.
\end{proof}

\smallskip

\begin{proof}{\bf{}Proof of Lemma~\ref{main-dist-est}.}
First we establish an $L^{2}$ bound for $|x|^{-1}d(x)$, namely
\[%
\int_{M\cap B_{R}}|x|^{-2}d^{2}(x,y)\,d\mu(x,y) \le C R^{-2}\int_{M\cap B_{R/2}\setminus\{
  (x,y):|x|<\smfr{1}{4}R\}}d^{2}(x,y)\,d\mu(x,y) %
\pdl{pf-m-1}
\]%
for $R$ sufficiently large (depending on $M$), where $C=C(\lambda,\Lambda,\C)$.

To check this observe that since $R^{-n-\ell-2}\int_{M\cap B_{R}}d^{2}\to 0$ by~\ref{tangent-cone}, we can apply
Lemma~\ref{tech-prelim} for sufficiently large $R$ giving conclusions~(i),(ii) of~\ref{tech-prelim} with
$\delta=\fr{1}{8}$ and $\theta=\fr{1}{2}$.

Let $\epsilon=\epsilon(\delta,\lambda,\Lambda,\theta)\in(0,\epsilon_{0}]$ be as in Lemma~\ref{tech-prelim} with
$\delta=\fr{1}{8}$ and $\theta=\fr{1}{2}$, and let $\bigl(d/|x|\bigr)_{\epsilon}=\min\{d(x)/|x|,\epsilon\}$. Then,
since $d(x)/|x|\le 1$ for all $x$,
\[%
d(x)/|x| \le  \epsilon^{-1}\bigl(d(x)/|x|\bigr)_{\epsilon} \text{ at all points of } M, %
\pdl{pf-m-2}
\]%
and since $\nabla_{M}\bigl(d(x)/|x|\bigr)_{\epsilon}=0$ for $\mu$-a.e.\ point $(x,y)$ with $d(x)/|x|\ge \epsilon$,
by~\ref{tech-prelim}(ii)
{\abovedisplayskip8pt\belowdisplayskip8pt%
\[%
|(x,0)\cdot\nabla_{M}\bigl(d(x)/|x|\bigr)_{\epsilon}| \le 2|x|^{-1}|(x,0)\cdot\nu(x,y)| %
                                                                      \text{ at $\mu$-a.e.\ point of } M\setminus B_{R_{0}}, %
\pdl{pf-m-3}
\]%
for suitable $R_{0}=R_{0}(M,\lambda,\Lambda,\C)$.

Take $R>R_{0}$ and $\zeta\in C^{1}_{c}(\breve B^{n+1}_{R})$, $\chi\in C^{1}_{c}(\breve B^{\ell}_{R})$ with
$\zeta(x)\equiv 1$ for $|x|\le \fr{1}{4} R$, $\zeta(x)=0$ for $|x|>\ha R$, $\chi(y)\equiv 1$ for $|y|\le \fr{1}{4}R$,
$\chi(y)=0$ for $|y|>\ha R$, and $|D\zeta|,|D\chi|\le 3 R^{-1}$.  Since
\[%
 \dvg_{M}(x,0)\ge n> 2\text{ and } |(x,0)\cdot\nabla_{M}\chi(y)|\le |x|\,|D\chi(y)|\,|\nu_{y}|, 
\]%
we can apply~\ref{pf-m-2}, \ref{pf-m-3}, and~\ref{stationarity} with
$Z_{|(x,y)}=\bigl(d(x)/|x|\bigr)_{\epsilon}^{2}\zeta^{2}(x)\chi^{2}(y)(x,0)$ to conclude
\begin{align*}%
&\avint_{M\cap B_{R/4}}\hskip-8pt |x|^{-2}d^{2}(x)\,d\mu \ptg{pf-m-4}\\ %
&\hskip.2in \le  C\avint_{M\cap B_{R}}\hskip-6pt\bigl(|x|^{-2}((x,0)\cdot \nu)^{2}+|\nu_{y}|^{2}\bigr)\,d\mu+ %
CR^{-2}\avint_{M\cap B_{R}\setminus \{(x,y):|x|<\frac{1}{4} R\}}\hskip-20pt d^{2}\,d\mu \\
&\hskip0.2in \le 2 C\avint_{M\cap B_{R}}\hskip-6pt \bigl(|x|^{-2}((x,y)\cdot
\nu)^{2}+R^{2}|x|^{-2}|\nu_{y}|^{2}\bigr)\,d\mu+ \\  
\noalign{\nobreak\vskip-4pt}
&\hskip2.6in CR^{-2}\avint_{M\cap B_{R}\setminus\{(x,y):|x|<\frac{1}{4}R\}}\hskip-20pt d^{2}\,d\mu,
\end{align*}%
where $C=C(\lambda,\Lambda,\C)$.

By~\ref{strong-doub-lem}, with either of the choices $w=(x,y)\cdot \nu$ or $w=R\nu_{y_{j}}$, 
\[%
\avint_{M\cap B_{R}}|x|^{-2}w^{2}\le CR^{-2}\avint_{M\cap B_{R/2}\setminus\{(x,y):|x|<\smfr{1}{3} R\}}w^{2}, %
\]%
and hence~\ref{pf-m-4} gives
\begin{align*}%
  &\avint_{M\cap B_{R/4}}\hskip-8pt |x|^{-2}d^{2}(x)\,d\mu\le %
  CR^{-2}\avint_{M\cap
B_{R}\setminus\{(x,y):|x|<\smfr{1}{3}R\}}\hskip-20pt(\,((x,y)\cdot\nu)^{2}+R^{2}|\nu_{y}|^{2}\,)\ptg{pf-m-4.1}\\ %
  &\hskip2.5in +CR^{-2}\avint_{M\cap B_{R}\setminus\{(x,y):|x|<\smfr{1}{4}R\}}\hskip-30pt d^{2}\,d\mu. %
\end{align*}%
In view of~\ref{tangent-cone} there exists $R_{0}>1$ and $\delta,\tilde\delta:(0,\infty)\to (0,\infty)$ such that
$\delta(t),\tilde\delta(t)\to 0$ as $t\to\infty$ and 
\[%
\left\{\begin{aligned}%
& \bigl\{(x,y)\in M:|x|\le \tilde\delta(|(x,y))|(x,y)|\bigr\}\setminus B_{R_{0}}\subset \cup_{j=1}^{q}\graph u_{j} \subset M \\ %
&{\sup}_{(\xi,\eta)\in \Omega\setminus B_{R}}\tsum_{j=1}^{q}\bigl(|Du_{j}(\xi,\eta)|+|(\xi,\eta)|^{-1}u_{j}(\xi,\eta)\bigr) \to 0%
 \text{ as }R\to\infty,  %
\end{aligned}\right.%
\pdl{pf-m-10}
\]%
where $u_{j}$ are positive $C^{2}$ functions on the domain $\Omega$,
\[%
\Omega\supset\{(x,y)\in\C:|x|<\delta((|x,y|))|(x,y)|\}\setminus B_{R_{0}}.
\]}%
For $(x,y)=(\xi+u_{j}(\xi,y)\nu_{\C_{0}}(\xi),y)$ with $(\xi,y)\in\Omega$, take an orthonormal basis
$\tau_{1},\ldots,\tau_{n+\ell}$ for $T_{(\xi,y)}\C$ with $\tau_{1},\ldots,\tau_{n-1}$ principal directions for $\C_{0}\cap
\Sph^{n}$ (so $\nabla_{\tau_{i}}\nu_{\C_{0}}=\kappa_{i}\tau_{i}$ for $i=1,\ldots,n-1$), $\tau_{n}=|\xi|^{-1}\xi$, and
$\tau_{n+j}=e_{n+1+j}$, $j=1,\ldots,\ell$.  Then the unit normal $\nu(x,y)$ of $M$ is
{\abovedisplayskip8pt\belowdisplayskip8pt%
\begin{align*}%
&\nu(x,y)= \bigl(1+\tsum_{i=1}^{n}(1+\kappa_{i}u_{j}(\xi,y))^{-2}(D_{\tau_{i}}u_{j}(\xi,y))^{2}+ %
|D_{y}u(\xi,y)|^{2}\bigr)^{-1/2} \\ %
\noalign{\vskip-1pt}    %
&\hskip0.2in \times\bigl(\nu_{\C_{0}}(\xi,y)-\tsum_{i=1}^{n}(1+\kappa_{i}u_{j}(\xi,y))^{-1} %
D_{\tau_{i}}u_{j}(\xi,y)\tau_{i}  -\tsum_{k=1}^{\ell}D_{y_{k}}u_{j}(\xi,y)e_{n+1+k}\bigr), %
\end{align*}%
so for $R$ sufficiently large
\[%
|\nu_{y}(x,y)|\le |D_{y}u_{j}(\xi,y)|\text{ and  } |(x,0)\cdot\nu_{M}(x,y)|\le |u_{j}(\xi,y)|+2|\xi||\nabla_{\C}u_{j}(\xi,y)|. %
\pdl{nu-ineq}
\]%
Also, since $u_{j}$ satisfies the equation
$\mathcal{M}_{\C}(u_{j})=0$, where $\mathcal{M}_{\C}$ satisfies~\ref{L-M-3} with $\C$ in place of $M$, we have
\[%
\int_{B_{R}\setminus\{(x,y):|x|<\smfr{1}{3} R\}}|\nabla_{\C}u_{j}|^{2}\le CR^{-2}
\int_{B_{3R/2}\setminus\{(x,y):|x|<\smfr{1}{4} R\}}u_{j}^{2},
\pdl{ell-est}
\]%
Also, $d(\xi+u_{j}(\xi,y)\nu_{\C_{0}}(\xi))=u_{j}(\xi,y)$ so, by~\ref{nu-ineq} and~\ref{ell-est}, 
\[%
\avint_{M\cap B_{R}\setminus\{(x,y):|x|<\smfr{1}{3} R\}}w^{2}\le %
                           C\avint_{M\cap B_{3R/2}\setminus\{(x,y):|x|<\smfr{1}{4}R\}}d^{2} %
\pdl{pf-m-7}
\]%
for either of the choices $w=(x,y)\cdot \nu$ and $w=R\nu_{y_{j}}$, and hence~\ref{pf-m-4.1} implies
\[%
\avint_{M\cap B_{R/4}}|x|^{-2}d^{2}(x)\le CR^{-2}\avint_{M\cap B_{3R/2}\setminus\{(x,y):|x|<\smfr{1}{4}R\}}d^{2},
\quad C=C(\lambda,\Lambda,\C).
\pdl{pf-m-7.1}
\]%
Since $\mathcal{M}(u_{j})=0$ and the $u_{j}$ are positive with small gradient, and also $d(\xi +
u_{j}(\xi,\eta)\nu_{\C_{0}}(\xi),\eta)=u_{j}(\xi,\eta)$, we can use the Harnack inequality in balls of radius $R/20$, to give
\[%
\avint_{M\cap B_{3R/2}\setminus\{(x,y):|x|<\smfr{1}{4}R\}}d^{2}\le
C\avint_{M\cap B_{R/8}\setminus\{(x,y):|x|<\smfr{1}{16}R\}}d^{2},
\]%
and so~\ref{pf-m-7.1} gives 
\[%
\avint_{M\cap B_{R/4}}|x|^{-2}d^{2}(x) \le CR^{-2}\avint_{M\cap B_{R/8}\setminus\{(x,y):|x|<\smfr{1}{16}R\}}d^{2},
\quad C=C(\lambda,\Lambda,\C).
\pdl{pf-m-7.2}
\]%
Then~\ref{pf-m-1} follows from~\ref{pf-m-7.2} after replacing $R$ by $4R$. \ref{pf-m-1} in particular implies, for all
sufficiently large $R$,
\[%
\left\{\begin{aligned}%
&\int_{M\cap B_{R}\cap \{(x,y):|x|<\delta R\}}d^{2}\,d\mu\le C\delta^{2}\int_{M\cap B_{R/2}} d^{2}\,d\mu %
\,\,\,\forall  \delta\in (0,\ha], \\ %
&\int_{M\cap B_{R}}d^{2}\,d\mu\le C\int_{M\cap B_{R/2}\setminus\{(x,y):|x|<\frac{1}{4}R\}} d^{2}\,d\mu, %
\end{aligned}\right.%
\pdl{pf-m-7a}
\]%
where $C=C(\lambda,\Lambda,\C)$.  

So let $R_{k}\to\infty$ be arbitrary, $M_{k}=R_{k}^{-1}M$,\,$d_{k}=d|M_{k}$, and
\[%
u^{k}_{j}(x,y)=R_{k}^{-1}u_{j}(R_{k}x,R_{k}y)\big/\bigl(\smavint_{M_{k}\cap B_{1}}d_{k}^{2}\bigr)^{1/2}.
\]}%
By virtue of~\ref{pf-m-1}, \ref{pf-m-10}, and \ref{pf-m-7a}, we have a subsequence of $k$ (still denoted $k$) such that
$u^{k}_{j}$ converges locally in $\C$ to a solution $v_{j}\ge 0$ of $\mathcal{L}_{\C}v_{j}=0$ with $v_{j}$ strictly positive
for at least one $j$, and $\int_{\C\cap B_{R}}|x|^{-2}v_{j}^{2}\,d\mu<\infty$ for each $R>0$.  Hence $v_{j}$ has a
representation of the form~\ref{LC-soln} on all of $\C$, and then since $v_{j}\ge 0$ we must have
$\smash{v_{j}=c_{j}r^{\gamma_{1}^{+}}\varphi_{1}}$ with $c_{j}\ge 0$ and $c_{j}>0$ for at least one $j$.  But then
$\smavint_{\C\cap B_{R}}\sum_{j=1}^{q}v^{2}_{j}$ is constant, independent of $R$, and so (using~\ref{pf-m-7a} again)
\[%
\lim_{k\to\infty}\avint_{M_{k}\cap B_{2R_{k}}}d^{2}\bigg/\avint_{M_{k}\cap B_{R_{k}}}d^{2}=
\avint_{\C\cap B_{2}}\sum_{j=1}^{q}v^{2}_{j}\bigg/\avint_{\C\cap B_{1}}\sum_{j=1}^{q}v^{2}_{j}= 1.
\]%
Therefore, in view of the arbitrariness of the sequence $R_{k}$, for any $\alpha\in (0,\ha]$ we have
\[%
2^{-\alpha}\le \avint_{M\cap B_{R}}d^{2}\bigg/\avint_{M\cap B_{R/2}}d^{2}   \le 2^{\alpha}\,\,\, \forall\,R>R_{0}, 
\]%
with $R_{0}=R_{0}(M,\lambda,\Lambda,\alpha)$, and~\ref{main-dist-est} follows. \end{proof}

\section{Proof of Theorem~\ref{main-th}} \label{th-1-pf}

\noindent According to Lemma~\ref{tangent-cone} the unique tangent cone of
$M$ at $\infty$ is $\C$ (with multiplicity $q$).  Let $\alpha\in (0,1)$.  By Lemma~\ref{main-dist-est}
\[%
\avint_{M\cap B_{R}}d^{2}\,d\mu \le  R^{\alpha} %
\dl{pf-th-1}
\]%
for all sufficiently large $R$.  Also by~\ref{growth-lem} either $\nu_{y}$ is identically zero or there is a lower bound
\[%
\avint_{M\cap B_{R}}|\nu_{y}|^{2}\ge R^{-\alpha} %
\dl{pf-th-2}
\]%
for all sufficiently large $R$, and by Lemma~\ref{strong-doub-lem} with $w=\nu_{y_{j}}$ we have
\[%
\avint_{M\cap B_{R}}|\nu_{y}|^{2}\le C\avint_{M\cap B_{R/2}\setminus\{(x,y):|x|<\smfr{1}{3} R\}}|\nu_{y}|^{2}, %
\quad C =C(q,\lambda,\C). %
\dl{pf-th-3}
\]%
By inequality~\ref{pf-m-7} in the proof of Lemma~\ref{main-dist-est} (with $R/2$ in place of $R$), 
\[%
\avint_{M\cap B_{R/2}\setminus\{(x,y):|x|<\smfr{1}{3} R\}}|\nu_{y}|^{2}\le CR^{-2}\avint_{M\cap B_{R}}d^{2}. %
\dl{pf-th-3.1}
\]%
Combining~\ref{pf-th-1}, \ref{pf-th-2}, \ref{pf-th-3} and~\ref{pf-th-3.1}, we then have
\[%
R^{-\alpha}\le \avint_{M\cap B_{R}}|\nu_{y}|^{2} \le CR^{-2} \avint_{M\cap B_{R}}d^{2} \le  R^{-2+\alpha}%
\]%
for all sufficiently large $R$ (depending on $\alpha,q,\lambda,\C$), which is impossible for $R>1$. Thus the alternative
that $\nu_{y}$ is identically zero on $M$ must hold, and so $M$ is a cylinder $S\times\R^{\ell}$.\nobreak

\section{Proof of Corollary~\ref{co-1}} \label{proof-liou}

\noindent We aim here to show that the hypotheses of Corollary~\ref{co-1} ensure that $M$ is strictly stable as in~\ref{str-stab-M},
so that Corollary~\ref{co-1} is then implied by Theorem~\ref{main-th}.

Before discussing the proof that $M$ is strictly stable, we need a couple of preliminary results.  First recall that, according to
\cite[Theorem 2.1]{HarS85}, since $\C_{0}$ is minimizing there is a smooth embedded minimal hypersurface $S\subset U_{+}$
($U_{+}$ either one of the two connected components of $\R^{n+1}\setminus\overline{\C}_{0}$), with
\begin{align*}%
&\text{$\dist(S,\{0\})=1$, $\sing S=\emptyset$ and $S$ is a smooth \emph{radial graph}; i.e.\ each }\dtg{props-S} \\
\noalign{\vskip-2pt}
&\hskip0.4in\text{ray  $\bigl\{tx:t>0\bigr\}$ with $x\in U_{+}$ meets $S$ transversely in just one point, }
\end{align*}%
so in particular
\begin{align*}%
 &x\cdot \nu_{S}(x)>0 \text{ for each $x\in S$, where $\nu_{S}$ is the smooth} \dtg{radial-graph}\\ %
  \noalign{\nobreak\vskip-2pt}
 &\hskip1.2in   \text{ unit normal of $S$ pointing away from $\C_{0}$.} %
\end{align*}%
Furthermore, since $\C_{0}$ is strictly stable and strictly minimizing, \cite[Theorem 3.2]{HarS85} ensures that, for
$R_{0}=R_{0}(\C_{0})$ large enough, there is an $C^{2}$ function $u$ defined on an open subset of
$\Omega\subset\C_{0}$ with $\Omega\supset \bigl\{x\in\C_{0}:|x|>2R_{0}\bigr\}$ and
\[%
S\setminus B^{n+1}_{R_{0}}=\bigl\{x+u(x)\nu_{\C_{0}}(x):x\in\Omega\bigr\} \text{ with }u(x)=  %
\kappa|x|^{\gamma_{1}}\varphi_{1}(|x|^{-1}x) +E(x),  %
\dl{S-asymp}
\]%
where $\kappa$ is a positive constant, $\nu_{\C_{0}}$ is the unit normal of $\C_{0}$ pointing into $U_{+}$,
$\varphi_{1}>0$ is as in~\ref{char-exps} with $j=1$, and, for some $\alpha=\alpha(\C_{0})>0$,
\[%
\lim_{R\to\infty}\sup_{|x|>R}|x|^{k+|\gamma_{1}|+\alpha}|D^{k}E(x)|= 0\text{ for }k=0,1,2.
\]%

We claim that $S$ is strictly stable:

\begin{state}{\bf{}\tl{strict-stab-S} Lemma.}%
If $S$ is as above, there is a constant $\lambda=\lambda(\C_{0})>0$ such that
\[%
\lambda\int_{S}|x|^{-2}\zeta^{2}(x,y)\,d\mu(x,y)\le  %
                        \int_{S}\bigl(\bigl|\nabla_{S}\zeta\bigr|^{2}-|A_{S}|^{2}\zeta^{2}\bigr)\,d\mu, %
\quad \zeta\in C_{c}^{1}(\R^{n+1}), %
\]%
where $|A_{S}|$ is the length of the second fundamental form of $S$.
\end{state}%

\begin{proof}{\bf{}Proof:}
The normal part of the velocity vector of the homotheties $\lambda S|_{\lambda>0}$ at $\lambda=1$ is
$\psi=x\cdot\nu_{S}(x)\, (\,>0 \text{ by~\ref{radial-graph}})$, and since the $\lambda S$ are minimal hypersurfaces, this is a
Jacobi function, i.e.\ a solution of
\[%
\Delta_{S}\psi+|A_{S}|^{2}\psi=0. 
\pdl{jac-for-S}
\]%
By properties~\ref{radial-graph} and ~\ref{S-asymp} we also have
\[%
C^{-1}|x|^{\gamma_{1}}\le \psi(x)\le C|x|^{\gamma_{1}} \text{ on $S$},\quad C=C(\C_{0}).
\pdl{bounds-for-psi}
\]%
After a direct calculation and an integration by parts,~\ref{jac-for-S} implies
\[%
\int_{S}\psi^{2} |\nabla_{\!S}\bigl(\zeta/\psi\bigr)|^{2}\,d\mu = %
                      \int_{S}\bigl(|\nabla_{\!S} \zeta|^{2}-|A_{S}|^{2}\zeta^{2}\bigr)\,d\mu,  %
\pdl{quot}
\]%
and using the first variation formula~\ref{stationarity} with $Z(x)=|x|^{-p-2}f^{2}(x)x$ and noting that $\dvg_{S}x=n$,
we have, after an application of the Cauchy-Schwarz inequality,
\[%
\int_{S} |x|^{-p-2}f^{2} \le C\int_{S} |x|^{-p}|\nabla_{S}f|^{2},\,\, f\in C^{1}_{c}(\R^{n+1}),\quad p<n-2, \,\, C=C(p,n). %
\pdl{1st-var}
 \]%
 \ref{strict-stab-S} is now proved by taking $f=\zeta/\psi$ and $p=2|\gamma_{1}|\,(<n-2)$, and using~\ref{bounds-for-psi},
\ref{quot} and~\ref{1st-var}.
\end{proof}

\smallskip

As a corollary,  any hypersurface sufficiently close to $S$ in the appropriately scaled $C^{2}$ sense must also be strictly
stable:

\begin{state}{\bf{}\tl{co-stab-S} Corollary.}%
For each $\theta\in(0,\ha]$, there is $\delta=\delta(\C_{0},\theta)>0$ such that if $v\in C^{2}(S)$,
$M_{0}=\{x+v(x)\nu_{S}(x):x\in S\} \text{ and } |x||\nabla^{2}_{S}v|+|\nabla_{S}v|+|x|^{-1}|v|\le \delta\,\,\forall x\in S$, then
$M_{0}$ satisfies the inequality
\[%
\lambda(1-\theta)\int_{M_{0}}|x|^{-2}\zeta^{2}\,d\mu\le
\int_{M_{0}}\bigl(\bigl|\nabla_{M_{0}}\zeta\bigr|^{2}-|A_{M_{0}}|^{2}\zeta^{2}\bigr)\,d\mu, %
\quad \zeta\in C_{c}^{1}(\R^{n+1}),
\]%
where $|A_{M_{0}}|$ is the length of the second fundamental form of $M_{0}$, and $\lambda$ is the constant
of~{\rm\ref{strict-stab-S}}. 
\end{state}%

\begin{proof}{\bf{}Proof:}  By~\ref{strict-stab-S}, with $\tilde\zeta(x)=\zeta(x+v(x)\nu_{S}(x))$ for $x\in S$,
\[%
\lambda \int_{S} |x|^{-2}\tilde\zeta^{2}\,d\mu \le %
\int_{S}\bigl(\bigl|\nabla_{S}\tilde\zeta\bigr|^{2}-|A_{S}|^{2}\tilde\zeta^{2}\bigr)\,d\mu %
\pdl{co-stab-1}
\]%
and for any $C^{1}$ function $f$ with compact support on $M_{0}$, with $\tilde f(x)=f(x+v(x)\nu_{S}(x))$ for $x\in S$,
\[%
\left\{\hskip2pt\begin{aligned}%
&\tint_{S}\tilde f\,d\mu = \tint_{M_{0}} f J\,d\mu  \text{ with } |J-1|\le C\delta \text{ (by the area
formula)},\\ & |\nabla_{S}\tilde f(x)-(\nabla_{M_{0}}f)(x+v(x)\nu_{S}(x)|\le C\delta
|\nabla_{S}\tilde f(x)|\\ & |A_{S}(x)-A_{M_{0}}(x+v(x)\nu_{S}(x))|< C|x|^{-1}\delta,\,\,\,  C^{-1}|x|^{-1}\le
|A_{S}(x)|\le C|x|^{-1},
\end{aligned}\right.%
\pdl{co-stab-2}
\]%
where $C=C(S)$. By combining ~\ref{co-stab-1} and~\ref{co-stab-2} we then have the required inequality with
$\theta=C\delta$.  \end{proof}

\smallskip

Next we need a uniqueness theorem for stationary integer multiplicity varifolds with support contained in
$\overline U_{+}$.

\begin{state}{\bf{}\tl{pre-Liouville} Lemma.}%
If $M_{0}$ is a stationary integer multiplicity $n$-dimensional varifold in $\R^{n+1}$ with the properties that
$\spt M_{0}\subset\overline U_{+}$, $\spt M_{0}\neq \overline\C_{0}$, and $ \sup_{R>1}R^{-n}\mu(M_{0}\cap
B_{R})<2\mu(\C_{0}\cap B_{1})$, then $M_{0}=\lambda S$ (with multiplicity~$1$) for some $\lambda>0$, where $S$ is
as in~{\rm\ref{props-S}}.
\end{state}%

\begin{proof}{\bf{}Proof of \ref{pre-Liouville}:}
Let $C(M_{0})$ be a tangent cone of $M_{0}$ at $\infty$. Then by the Allard compactness theorem $C(M_{0})$ is a stationary
integer multiplicity varifold with $\mu_{C(M_{0})}(B_{1})<2\mu(\C_{0}\cap B_{1})$ and $\spt C(M_{0})\subset \overline
U_{+}$. If $\spt C(M_{0})$ contains a ray of $\C_{0}$ then, by the Solomon-White maximum principle~\cite{SolW89}, either
$C(M_{0})=\C_{0}$ (with multiplicity one) or else $C(M_{0})=\C_{0}+V_{1}$, where $V_{1}$ is a non-zero integer multiplicity
cone with $\mu(V_{1}\cap B_{1}(0))<\mu(\C_{0}\cap B_{1})$ and support contained in $\overline U_{+}$. On the other hand if
$\spt C(M_{0})\cap \C_{0}=\emptyset$ then a rotation of $\spt C(M_{0})$ has a ray in common with $\C_{0}$ so by the same
argument applied to this rotation we still conclude that there is a stationary cone $V_{1}$ with
$\mu_{V_{1}}(B_{1})<\mu(\C_{0}\cap B_{1})$ and $\spt V_{1}\subset \overline U_{+}$.  But now by applying exactly the same
reasoning to $V_{1}$ we infer $\mu_{V_{1}}(B_{1})\ge \mu(\C_{0}\cap B_{1})$, a contradiction. So $C(M_{0})=\C_{0}$ and
hence $\C_{0}$, with multiplicity one, is the unique tangent cone of $M_{0}$ at infinity.  Hence there is a $R_{0}>1$ and a
$C^{2}(\C_{0}\setminus B_{R_{0}})$ function $h$ with
\begin{align*}%
&\hskip0.7in\text{$\sup_{x\in \C_{0}\setminus B_{R}}(|Dh(x)|+|x|^{-1}h(x))\to 0$ as $R\to \infty$}, \ptg{def-h-1}\\ %
\noalign{\vskip3pt}
&\{x+h(x)\nu_{\C_{0}}(x):x\in \C_{0}\setminus B_{R_{0}}\}\subset M_{0} \ptg{def-h-2} \\ %
\noalign{\vskip-2pt} %
&\hskip1in \text{ and }M_{0}\setminus \{x+h(x)\nu_{\C_{0}}(x):x\in \C_{0}\setminus B_{R_{0}}\} \text{ is compact.} %
\end{align*}%
We also claim 
\[%
0\notin M_{0}. %
\pdl{0-notin-M}
\]%
Indeed otherwise the Ilmanen maximum principle~\cite{Ilm96} would give $M_{0}\cap \C_{0}\neq \emptyset$ and the
above argument using the Solomon-White maximum principle can be repeated, giving $M_{0}=\C_{0}$, contrary to the
hypothesis that $\spt M_{0}\neq \overline\C_{0}$.

Observe next that if $r_{k}\to\infty$ the scaled minimal hypersurfaces $r_{k}^{-1}M_{0}$ are represented by the graphs
(taken off $\C_{0}$) of the functions $r_{k}^{-1}h(r_{k}x)\to 0$, and hence, for any given $\omega_{0}\in\Sigma\,
(=\C_{0}\cap\partial B_{1})$, the rescalings $(h(r_{k}\omega_{0}))^{-1}h(r_{k}r\omega)$ (which are bounded above and
below by positive constants on any given compact subset of $\C_{0}$ by~\ref{def-h-1} and the Harnack inequality) generate
positive solutions of the Jacobi equation $\mathcal{L}_{\C_{0}}v=0$ on $\C_{0}$ as $k\to\infty$.  But,
by~\ref{eigenvals} and \ref{exp-gamma}, $(c_{1}r^{\gamma_{1}}+c_{2}r^{\gamma_{1}^{-}})\,\varphi_{1}(\omega)$ with
$c_{1},c_{2}\ge 0$ and $c_{1}+c_{2}>0$ are the only positive solutions of $\mathcal{L}_{\C_{0}}(\varphi)=0$ on all of
$\C_{0}$.  Thus in view of the arbitrariness of the sequence $r_{k}$, we have shown that
\[%
h(r\omega) =c(r) \varphi_{1}(\omega) +o(c(r)) \text{ as $r\to \infty$, uniformly for $\omega\in\Sigma$}, %
\pdl{M-as}
\]%
and hence there are 
\[%
c_{-}(r)<c(r)<c_{+}(r)\text{ with }c_{-}(r)\varphi_{1}(\omega)<h(r\omega)<c_{+}(r)\varphi_{1}(\omega)  %
\text{ and }c_{+}(r)/c_{-}(r)\to 1.  %
\]%
Now, for suitable $R_{0}>0$, $S\setminus B_{R_{0}}$ ($S$ as in~\ref{props-S}) also has a representation of this form with
some $\tilde h$ in place of $h$, where
$$%
\tilde h(r\omega)=\kappa r^{\gamma_{1}}\varphi_{1}(\omega) + o(r^{\gamma_{1}}) \text{ as $r\to \infty$,
uniformly in $\omega$}
$$%
and, similar to the choice of $c_{\pm}$, we can take  $\tilde c_{\pm}(r)$ such that 
\[%
\tilde c_{-}(r)\varphi_{1}(\omega)<\tilde h(r\omega)< %
\tilde c_{+}(r)\varphi_{1}(\omega)\text{ and }\tilde c_{+}(r)/\tilde c_{-}(r)\to 1\text{ as }r\to \infty.  %
\]%
Now $\lambda S\setminus B_{\lambda R_{0}}$ can be represented by the geometrically scaled function $\tilde h_{\lambda}$
with 
\[%
\tilde h_{\lambda}(r\omega)=\kappa \lambda^{1+|\gamma_{1}|}r^{-|\gamma_{1}|}\varphi_{1} + 
o(\kappa \lambda^{1+|\gamma_{1}|}r^{-|\gamma_{1}|}\varphi_{1})
\pdl{M-as-lam}
\]%
and we let $\lambda_{k}^{-}$ be the largest value of $\lambda$ such that $\tilde h_{\lambda}(r_{k}\omega)\le
c_{-}(r_{k})\varphi_{1}(\omega)$ and $\lambda_{k}^{+}$ the smallest value of $\lambda$ such that $\tilde
h_{\lambda}(r_{k}\omega)\ge c_{+}(r_{k})\varphi_{1}(\omega)$. Evidently there are then points $\omega_{\pm} \in
\Sigma$ with $\tilde h_{\lambda_{k}^{\pm}}(r_{k}\omega_{\pm})=c_{\pm}(r_{k})\varphi_{1}(\omega_{\pm})$ respectively. 
Also $M_{0}\cap \breve B_{r_{k}}$ must entirely lie in the component $B_{r_{k}}\setminus \lambda_{k}^{+} S$ which
contains $\C_{0}\cap B_{r_{k}}$; otherwise we could take the smallest $\lambda>\lambda_{k}^{+}$ such that $M_{0} \cap
B_{r_{k}}$ lies on in the closure of that component of $B_{r_{k}}\setminus \lambda S$, and $M_{0}\cap \breve
B_{r_{k}}\cap \lambda S\neq \emptyset$, which contradicts the maximum principle~\cite{SolW89}.  Similarly, $M_{0}\cap
\breve B_{r_{k}}$ must entirely lie in the component $B_{r_{k}}\setminus \lambda_{k}^{-} S$ which does not contain
$\C_{0}\cap B_{r_{k}}$.

Thus $M_{0}\cap B_{r_{k}}$ lies between $\lambda_{k}^{+}S$ and $\lambda_{k}^{-}S$ and by construction
$\lambda_{k}^{+}/\lambda_{k}^{-}\to 1$, and since $\lambda_{k}^{-}$ is bounded above, we then have a subsequence of $k$
(still denoted $k$) such that $\lambda_{k}^{\pm}$ have a common (positive) limit $\lambda$.  So $\spt M_{0}\subset
\lambda S$ and hence $M_{0}=\lambda S$ with multiplicity~$1$ by the constancy theorem and the fact that $
\sup_{R>1}R^{-n}\mu(M_{0}\cap B_{R})<2\mu(\C_{0}\cap B_{1})$.
\end{proof}

\medskip

Finally we need to show that any $M$ satisfying the hypotheses of Corollary~\ref{co-1} with sufficiently small $\epsilon_{0}$
must be strictly stable; then (as mentioned at the beginning of this section) Corollary~\ref{co-1} follows from
Theorem~\ref{main-th}.

\smallskip

\begin{state}{\bf{}\tl{stab-lem} Lemma.}%
For each $\alpha,\theta\in (0,1)$, there is $\epsilon_{0}=\epsilon_{0}(\C,\alpha,\theta)\in (0,\ha]$ such that if $M\subset
U_{+}\times\R^{\ell}$, ${\sup}_{R>1}R^{-n-\ell}\mu(M\cap B_{R})\le (2-\alpha)\mu(\C\cap B_{1})$, and
${\sup}_{M}|\nu_{y}| < \epsilon_{0}$, then $M$ is strictly stable in the sense that
\[%
(1-\theta)\lambda\int_{M}|x|^{-2}\zeta^{2}\,d\mu\le %
\int_{M}\bigl(\bigl|\nabla_{M}\zeta\bigr|^{2}-|A_{M}|^{2}\zeta^{2}\bigr)\,d\mu,\quad\zeta\in C_{c}^{1}(\R^{n+1+\ell}), %
\]%
with $\lambda=\lambda(\C_{0})>0$ as in~{\rm\ref{strict-stab-S}}.
\end{state}%

\begin{proof}{\bf{}Proof:}
For $y\in\R^{\ell}$, define
{\abovedisplayskip8pt\belowdisplayskip8pt%
\[%
M_{y}= \lambda_{y}^{-1}(M-(0,y)),\quad  \lambda_{y}=\dist(M-(0,y),0).
\pdl{pf-0}
\]%
We claim that for each given $\delta>0$ the hypotheses of the lemma guarantee that a strip of $M$ near a given slice
$M_{z}=\{(x,y)\in M:y=z\}$ can be scaled so that it is $C^{2}$ close to $S\times B^{\ell}_{1}$ ($S$ as in~\ref{radial-graph}) in
the appropriately sense; more precisely, we claim that for each $\delta>0$ there is
$\epsilon_{0}=\epsilon_{0}(\C,\alpha,\theta,\delta)>0$ such that the hypotheses of the lemma imply that for each
$z\in\R^{\ell}$ there is $v_{z}\in C^{2}(\{(x,y):x\in S,\,|y|<1\})$ such that
\[%
\left\{\begin{aligned}
&M_{z}\cap \{(x,y):|y|<1\}  %
                             =\bigl\{(x+v_{z}(x,y)\nu_{S}(x),y):x\in S,\,|y|<1\bigr\}\\  %
&\, |x||\nabla^{2}_{S\times\R^{\ell}}v_{z}(x,y)|+|\nabla_{S\times\R^{\ell}}v_{z}|(x,y)+|x|^{-1}|v_{z}(x,y)|\le \delta %
                                                          \,\,\forall x\in S,|y|<1.  %
\end{aligned}\right.%
\pdl{pf-2}
\]%
Otherwise this fails with $\epsilon_{0}=1/k$ for each $k=1,2,\ldots$, so there are minimal submanifolds $M_{k}$ such
that the hypotheses hold with $\epsilon _{0}=1/k$ and with $M_{k}$ in place of $M$, yet there are points
$z_{k}\in M_{k}$ such that~\ref{pf-2} fails with $z_{k},M_{k}$ in place of $z,M$.

\vskip1pt

We claim first that then there are fixed $k_{0}=k_{0}(\delta)\in \{1,2,\ldots\}, \,R_{0}=R_{0}(\delta)>1$ such that
\[%
\dist((x,y),M_{z_{k}}) <\delta|x| \,\,\forall (x,y)\in \C\setminus B_{R_{0}}^{n+1}\text{ with }|x|\ge |y|,\,\,\,k\ge k_{0}. %
\pdl{pf-2z}
\]%
Otherwise there would be a subsequence of $k$ (still denoted $k$) with
\[%
\text{$\dist((x_{k},y_{k}),M_{z_{k}})\ge \delta |x_{k}|$ with $(x_{k},y_{k})\in\C$, $|y_{k}|\le |x_{k}|$ and $|x_{k}|\to
\infty$.}
\]%
Then let
\[%
\wtilde{M}_{k} = |x_{k}|^{-1}M_{z_{k}}.
\]%
By the Allard compactness theorem, $\wtilde{M}_{k}$ converges in the varifold sense to a stationary integer multiplicity
varifold $V$ with support $M$ and density $\Theta$, where $M$ is a closed rectifiable set, $\Theta$ is upper
semi-continuous on $\R^{n+1+\ell}$ and has integer values $\mu$-a.e.\ on $M$, with $\Theta=0$ on $\R^{n+1+\ell}\setminus
M$,  and
\[%
M=\{x:\Theta(x)\ge 1\},  %
\pdl{pf-2a}
\]%
\[%
\Theta(x) = \lim_{\rho\downarrow 0}\bigl(\omega_{n+\ell}\rho^{n+\ell}\bigr)^{-1}  %
                          \int_{M\cap B_{\rho}(x)}\Theta(\xi)\,d\mu(\xi)\,\,\,\forall x\in M, %
\pdl{pf-3}
\]%
where $\omega_{n+\ell}$ is the volume of the unit ball in $\R^{n+\ell}$, $0\in M\subset \overline
U_{+}\times\R^{\ell}$, and $M\cap U_{+}\times\R^{\ell}\neq 0$.  

Taking $x_{n+1+j}=y_{j}$, so points in $\R^{n+1+\ell}$ are written $x=(x_{1},\ldots,x_{n+1+\ell})$, and
letting $\nu_{k}=(\nu_{k\,1},\ldots,\nu_{k\,n+1+\ell})$ be a unit normal for $\wtilde{M}_{k}$ (so that the orthogonal
projection of $\R^{n+1+\ell}$ onto $T_{x}\wtilde{M}_{k}$ has matrix $(\delta_{ij}-\nu_{k\,i}\nu_{k\,j})$), the first
variation formula for $\wtilde{M}_{k}$ can be written
\[%
\int_{\wtilde{M}_{k}}\sum_{i,j=1}^{n+1+\ell}(\delta_{ij}-\nu_{k\,i}\nu_{k\,j})D_{i}X_{j}\,d\mu=0, %
                                          \,\,\, X_{j}\in C^{1}_{c}(\R^{n+1+\ell}),\,\,j=1,\ldots,n+1+\ell.%
\pdl{pf-4}
\]}%
Let $x_{0}\in M,\,\sigma>0$, and let $\tau\in C^{\infty}(\R^{n+1+\ell})$ with support $\tau\subset B_{\sigma}(x_{0})$,
and consider the function $T$ defined by
{\abovedisplayskip8pt\belowdisplayskip8pt%
\[%
T(x)=T(x_{1},\ldots,x_{n+1+\ell}) =  %
\int_{-\infty}^{x_{n+1+\ell}}\bigl(\tau(x_{1},\ldots,x_{n+\ell},t) %
                                      -\tau(x_{1},\ldots,2\sigma+x_{n+\ell},t)\bigr)\,dt. %
\]%
Evidently $T\in C^{\infty}_{c}(\R^{n+1+\ell})$ (indeed $T(x)=0$ on $\R^{n+1+\ell}\setminus K$, where $K$ is the cube
$\{x\in\R^{n+1+\ell}:|x_{i}-x_{0\,i}|<4\sigma\,\,\forall i=1,\ldots,n+1+\ell\}$).  Therefore we can use~\ref{pf-4} with
$X(x)=T(x)e_{n+1+\ell}$, giving
\[%
\int_{\wtilde{M}_{k}}\bigl(1-\nu_{k\,n+1+\ell}^{2}\bigr)\bigl( \tau(x)-\tau(x+2\sigma e_{n+1+\ell})\bigr)\,d\mu  %
    = \int_{\wtilde{M}_{k}}\sum_{i=1}^{n+\ell}\nu_{k\,n+1+\ell}\nu_{k\,i} D_{i}T\,d\mu, 
\]%
and hence 
\[%
\Bigl|\int_{\wtilde{M}_{k}} \tau(x)\,d\mu-\int_{\wtilde{M}_{k}}\tau(x+2\sigma e_{n+1+\ell})\,d\mu\Bigr| \le C/k. %
\]%
Using the fact that varifold convergence of $\wtilde{M}_{k}$ implies convergence of the corresponding mass measures
$\mu\res \wtilde{M}_{k}$ to $\mu\res \Theta$, we then have
\[%
\int_{M} \tau(x)\,\Theta(x) d\mu(x)=\int_{M}\tau(x+2\sigma e_{n+1+\ell})\,\Theta(x)d\mu(x). %
\pdl{pf-5}
\]%
Taking $\rho\in (0,\sigma)$ and replacing $\tau$ by a sequence $\tau_{k}$ with $\tau_{k}\downarrow
\chi_{B_{\rho}(x_{0})}$ (the indicator function of the closed ball $\smash{B_{\rho}(x_{0})}$), we then conclude 
\[%
\int_{M\cap B_{\rho}(x_{0})}\,\Theta\,d\mu=\int_{M\cap B_{\rho}(x_{0}+2\sigma e_{n+1+\ell})}\ \Theta\,d\mu  %
\]}%
and hence by~\ref{pf-3}
{\abovedisplayskip8pt\belowdisplayskip10pt%
\[%
\Theta(x_{0})=\Theta(x_{0}+2\sigma e_{n+1+\ell}). 
\]}%
In view of the arbitrariness of $\sigma$ this shows that $\Theta(x)$ is independent of the variable $x_{n+1+\ell}$, and the
same argument shows that $\Theta(x)$ is also independent of $x_{n+1+j}$, $j=1,\ldots,\ell-1$. Thus, by~\ref{pf-2a}, $M$ is
cylindrical: $M=M_{0}\times\R^{\ell}$, where $R^{-n}\int_{M_{0}\cap B_{R}}\Theta\,d\mu<2\mu(\C_{0}\cap B_{1})$,
$0\in M_{0}\subset \overline U_{+}$, and $M_{0}\cap U_{+}\neq \emptyset$.  Hence $M_{0}=\lambda S$ for some
$\lambda>0$ by virtue of Lemma~\ref{pre-Liouville}, contradicting the fact that $0\in M_{0}$. So~\ref{pf-2z} is proved,
and~\ref{pf-2z} together with the Allard regularity theorem implies that
$M_{z_{k}}\cap\bigl(\bigl\{(x,y)\in \R^{n+1+\ell}\setminus B_{2 R_{0}}:|y|<\ha|x| \bigr\}\bigr)$ 
 is $C^{2}$ close to $\C$,  in the sense that there is are $C^{2}$ functions $v_{k}$ on a domain in $\C$ with 
{\abovedisplayskip8pt\belowdisplayskip8pt%
\begin{align*}%
&\graph v_{k}=
M_{z_{k}}\cap\bigl(\bigl\{(x,y)\in \R^{n+1+\ell}\setminus B_{2 R_{0}}:|y|<\ha|x| \bigr\}\bigr)\ptg{pf-7}\\ %
&\hskip1.7in \text{ and }|x|^{-1}|v_{k}|+|\nabla_{\C}v_{k}|+|x||\nabla^{2}_{\C}v_{k}|<C\delta. 
\end{align*}}%

Next, exactly the same compactness discussion can be applied with $M_{z_{k}}$ in place of
$\wtilde{M}_{k}$, giving a cylindrical varifold limit $M=M_{0}\times\R^{\ell}$, $M_{0}$ a closed subset of $\overline
U_{+}$ with density $\Theta\ge 1$, but this time $\dist(M_{0},0)=1$.  Hence $\Theta\equiv 1$ and $M_{0}= S$
by~\ref{pre-Liouville}.  But then the Allard regularity theorem guarantees that the convergence of $M_{z_{k}}$ to $M$ is
smooth and hence, by virtue of~\ref{pf-7},~\ref{pf-2} holds with $M_{z_{k}}$ in place of $M$, a contradiction. Thus
\ref{pf-2} is proved.

Then by Corollary~\ref{co-stab-S}, for small enough $\delta=\delta(\C,\theta)$, for each $y\in\R^{\ell}$  
\[%
(1-\theta/2)\lambda\int_{M_{y}}|x|^{-2}\zeta_{y}^{2}(x)\,d\mu(x)\le
\int_{M_{y}}\bigl(\bigl|\nabla_{M_{y}}\zeta_{y}(x)\bigr|^{2}-|A_{M_{y}}|^{2}\zeta_{y}^{2}(x)\bigr)\,d\mu(x), %
\]%
$\zeta\in C_{c}^{1}(\R^{n+1})$, where $M_{y}=\{x\in\R^{ n+1}:(x,y)\in M\}$, $\zeta_{y}(x)=\zeta(x,y)$, and
$|A_{M_{y}}|$ is the length of the second fundamental form of $M_{y}$.

Since
\[%
||A_{M_{y}}|^{2}(x)-|A_{M}|^{2}(x,y)| \le C\delta/|x|^{2},\,\,\, |A_{M_{y}}|^{2}(x)\le C/|x|^{2}
\]%
by~\ref{pf-2}, the proof is completed by taking $\delta$ small enough, integrating with respect to $y\in \R^{\ell}$,
and using the coarea formula together with~\ref{pf-2}.
\end{proof}

\section{A Boundary Version of Theorem~\ref{main-th} and
Corollary~\ref{co-1}}\label{boundary-version}

\noindent Since it will be needed in~\cite{Sim21b}, we here also want to discuss a version of Theorem~\ref{main-th} which is
valid in case $M$ has a boundary.

\begin{state}{\bf{}\tl{bdry-ver} Theorem.}%
  Suppose $M\subset\R^{n+1}\times\{y:y_{\ell}\ge 0\}$ is a complete embedded minimal hypersurface-with-boundary, with
\[%
\left\{\hskip2pt\begin{aligned}%
&\partial M=S\times\{y:y_{\ell}=0\}, \\
\noalign{\vskip-1pt} %
&|\nu_{y_{\ell}}|<1\text{ on }S,\,\,  \\
\noalign{\vskip-1pt} %
&{\sup}_{R>1}R^{-n-\ell}\mu(M\cap B_{R}) <2\mu(\C\cap \bigl\{(x,y)\in B_{1}:y_{\ell}>0\bigr\}), \text{ and }\\ %
\noalign{\vskip-1pt} %
&M \subset U_{\lambda}\times\{y:y_{\ell}\ge 0\} %
\end{aligned}\right.%
\leqno{(\ddag)}
\]%
for some $\lambda\ge 1$, where $U_{\lambda}$ denotes the component of \,$U_{+}\setminus \lambda S$ with $\partial
U_{\lambda}=\overline\C_{0}\cup \lambda S$.

Suppose also that $M$ is strictly stable in the sense there is $\lambda>0$ such that that the inequality~\emph{\ref{str-stab-M}}
holds for all $\zeta\in C^{1}_{c}(M\cap \breve B_{R})$ with $e_{n+1+\ell}\cdot \nabla_{M}\zeta=0$ on $\partial M$ and for
all $R>0$.

Then 
\[%
M=S\times\{y:y_{\ell}\ge 0\}.
\]%
\end{state}%

\begin{proof}{\bf{}Proof:} Since $M\subset U_{\lambda}$, the Allard compactness theorem, applied to the rescalings $\tau
M,\,\tau\downarrow 0$, plus the constancy theorem, tells us that $M$ has a tangent cone at $\infty$ which is
$\C_{0}\times \{y:y_{\ell}\ge 0\}$ with some integer multiplicity $k\ge 1$, and then the condition $\sup_{R>1}R^{-n-1}\mu(M\cap
B_{R})<2\mu(\C\cap \bigl\{(x,y)\in B_{1}:y>0\bigr\})$ implies $k=1$. Thus $M$ has $\C_{0}\times\{y:y_{\ell}\ge 0\}$
with multiplicity~\!1 as its unique tangent cone at $\infty$.  

We claim that $\nu_{y_{\ell}}$ satisfies the free boundary condition
\[%
e_{n+1+\ell}\cdot \nabla_{M}\nu_{y_{\ell}}=0 \text{ at each point of $\partial M$. }%
\pdl{bdry-ver-1}
\]%
Indeed if $\Sigma\subset\R^{N}\times[0,\infty)$ is any minimal hypersurface-with-boundary with smooth unit normal
$\nu=(\nu_{1},\ldots,\nu_{N+1})$ and $|\nu_{N+1}|\neq 1$ on $\partial\Sigma$ (i.e.\ $\Sigma$ intersects the hyperplane
$x_{N+1}=0$ transversely), and if $\partial\Sigma$ is a minimal hypersurface in $\R^{N}\times\{0\}$, then
$e_{N+1}\cdot\nabla_{\Sigma} \nu_{N+1}=0$ on $\partial\Sigma$, as one readily checks by using the fact that the mean
curvature (i.e.\ trace of second fundamental form) of $\Sigma$ and $\partial\Sigma$ are both zero at each point of
$\partial\Sigma$.

We claim $\nu_{y_{\ell}}(x,y)=0\,\,\forall (x,y)\in M$.  To check this, first observe that there is a version
of~\ref{growth-lem} which is valid in the half-space $y_{\ell}\ge 0$ in the case when $w$ has free boundary condition
$e_{n+1+\ell}\cdot \nabla_{M}w=0$ on $\partial M$; indeed the proof of~\ref{growth-lem} goes through with little
change---the linear solutions $v$ of $\mathcal{L}_{\C}v=0$ obtained in the proof being defined on the half-space
$\C_{0}\times\{y:y_{\ell}\ge 0\}$ and having the free boundary condition $e_{n+1+\ell}\cdot\nabla_{\C}v=0$. So $v$
extends to a solution of $\mathcal{L}_{\C}v=0$ on all of $\C$ by even reflection and the rest of the argument is unchanged.

In particular, since $\nu_{y_{\ell}}$ has free boundary condition $0$ by~\ref{bdry-ver-1}, we have an analogue of
Lemma~\ref{strong-doub-lem} in the half-space, giving
\[%
R^{-\alpha}\le C\avint_{M\cap B^{+}_{R}}\nu_{y_{\ell}}^{2} %
\pdl{bdry-ver-2}
\]%
for each $\alpha\in(0,1)$, where $B^{+}_{R}=B_{R}\cap \{(x,y):y_{\ell}\ge 0\}$,
and, since $\nu_{y_{\ell}}\le 1$ there is a bound $\smavint_{M\cap B^{+}_{R}}\nu_{y_{\ell}}^{2}\le
CR^{n-2-\beta}$, so by Corollary~\ref{strong-doub-lem}
\[%
\avint_{M\cap B^{+}_{R}}\nu_{y_{\ell}}^{2} %
\le CR^{-\ell-2-\beta_{1}}\int_{B^{+}_{R/2}\setminus\{(x,y):|x|<\smfr{1}{3}R\}}\nu_{y_{\ell}}^{2}, %
\pdl{bdry-ver-3}
\]%
where $C=C(\C_{0},\alpha)$.  Since $M\subset U_{\lambda}\times[0,\infty)$ also we have
$d(x)\le C\lambda \min\{1,\,|x|^{\gamma_{1}}\}$  for all $x\in\R^{n}$ and hence
\[%
\avint_{M\cap B^{+}_{R}}d^{2} \le C\lambda^{2},  \quad C=C(\C_{0}).  %
\pdl{d-bound}
\]%
(Notice that here we do not need growth estimates for $\smavint_{B_{R}} d^{2}$ as in~\S\ref{dist-fn} because we are now assuming
$d\le C\lambda \min\{1,\,|x|^{\gamma_{1}}\}$.) Then the proof that $\nu_{y_{\ell}}$ must be identically zero is completed 
using~\ref{bdry-ver-2},  \ref{bdry-ver-3}, and \ref{d-bound} analogously to the proof of non-boundary version of
Theorem~\ref{main-th}. 

So $\nu_{y_{\ell}}$ is identically zero on $M$. This completes the proof in the case $\ell=1$ and shows that, for $\ell\ge 2$,
by even reflection $M$ extends to a minimal submanifold $\wtilde{M}$ (without boundary)
\[%
\wtilde{M}=M\cup\bigl\{(x,y_{1},\ldots,y_{\ell-1},-y_{\ell}):(x,y)\in M\bigr\}\subset U_{+}\times\R^{\ell}, 
\]%
which is translation invariant by translations in the direction of $e_{n+1+\ell}$. Then $\wtilde{M}$ is strictly stable and
Theorem~\ref{main-th} applies, giving $M=S\times \R^{\ell-1}\times[0,\infty)$, as claimed.\end{proof}

\smallskip

\begin{state}{\bf{}\tl{bdy-case} Corollary.}
There is $\delta=\delta(\C,\lambda)>0$ such that if $M$ satisfies~{\rm\ref{bdry-ver}}\hskip1pt$(\ddag)$ and
$\sup|\nu_{y}|<\delta$, then $M$ automatically satisfies the strict stability hypothesis in~{\rm\ref{bdry-ver}}, and hence
$M=S\times\{y:y_{\ell}\ge 0\}$.
\end{state}%

\begin{proof}{{\bf{}Proof:}}
With $\wtilde{M}_{k}$, $M_{z_{k}}$ as in the proof of~\ref{stab-lem}, with only minor modifications to the present situation
when $\partial M=S\times\{0\}$, we have $\wtilde M_{k}$ and $M_{z_{k}}$ both have cylindrical limits $S\times \R^{\ell}$
or $S\times \{y\in\R^{\ell}:y_{\ell}\ge K\}$ for suitable $K$, and so (using the Allard boundary regularity theorem in the latter
case) the $y=$ const.\ slices of $M$ have the $C^{2}$ approximation property~\ref{stab-lem}\ref{pf-2}. Hence, by integration
over the slices as in~\ref{stab-lem}, $M$ is strictly stable in the sense that~\ref{str-stab-M} holds for all $\zeta\in
C^{1}_{c}(\breve B_{R}^{+})$, where $\breve B_{R}^{+}=\{(x,y):|(x,y)|<R \text{ and }y_{\ell}\ge 0\}$.
\end{proof}

\section*{Appendix: Analyticity of $\beta$-harmonic functions}
\setcounter{sequation}{0}

\renewcommand{\thesequation}{A.\arabic{sequation}} %

\noindent For $\rho>0$, $r_{0}\ge 0$ and $y_{0}\in\R^{\ell}$, let 
\begin{align*}%
&B_{\rho}^{+}(r_{0},y_{0})= \{(r,y)\in \R\times\R^{\ell} :r\ge 0,\,\,(r-r_{0})^{2}+|y-y_{0}|^{2}\le \rho^{2}\}, \\ %
&\breve B^{+}_{\rho}(r_{0},y_{0})=\{(r,y)\in \R\times\R^{\ell}:r>0,\,\,(r-r_{0})^{2}+|y-y_{0}|^{2}<\rho^{2}\},
\end{align*}%
and $B^{+}_{\rho},\,\breve B^{+}_{\rho}$ will be used as abbreviations for $B^{+}_{\rho}(0,0),\, \breve B^{+}_{\rho}(0,0)$
respectively.

Our aim is to prove real analyticity extension across $r=0$ of $\beta$-harmonic functions, i.e.\ solutions $u\in
C^{\infty}(\breve B_{1}^{+})$ of
\[%
r^{-\gamma}\frac{\partial}{\partial r}\bigl(r^{\gamma}\frac{\partial u}{\partial r}\bigr) %
+\Delta_{y}u=0,
\dl{equn}
\]%
where $\gamma=1+\beta\,\,(\,>1)$ and $\Delta_{y}u=\tsum_{j=1}^{\ell}D_{y_{j}}D_{y_{j}}u$, assuming
\[%
\int_{B^{+}_{\rho}}(u_{r}^{2}+|u_{y}|^{2})\,\,r^{\gamma\!}drdy <\infty\,\,\, \forall \, \rho<1.%
\dl{w-1-2-bd}
\]%
In fact we show under these conditions that $u$ can be written as a convergent series of homogeneous $\beta$-harmonic
polynomials  in $\{(r,y):r\ge 0,\,r^{2}+|y|^{2}<1\}$ with the convergence uniform in $B_{\rho}^{+}$ for each $\rho<1$. 

Of course all solutions of~\ref{equn} satisfying \ref{w-1-2-bd} are automatically real-analytic in $\breve B^{+}_{1}$
because the equation is uniformly elliptic with real-analytic coefficients in each closed ball $\subset \breve B^{+}_{1}$. 
Also, if $\gamma$ is an integer $\ge 1$ then the operator in~\ref{equn} is just the Laplacian in $\R^{1+\gamma+\ell}$, at
least as it applies to functions $u=u(x,y)$ on $\R^{1+\gamma+\ell}$ which can be expressed as a function of $r=|x|$ and
$y$, so in this case smoothness across $r=0$ is Weyl's lemma and analyticity is then standard.

To handle the case of general $\gamma>1$, first note that by virtue of the calculus inequality
$\int_{0}^{a}r^{\gamma-2}f^{2}\,dr\le C\int_{0}^{a}r^{\gamma}(f'(r))^{2}\,dr$ for any $f\in C^{1}(0,a]$ with $f(a)=0$ (which is
proved simply by using the identity $\int_{0}^{a}\bigl(r^{\gamma-1}(\min\{|f|,k\})^{2}\bigr)'\,dr=0$ and then using
Cauchy-Schwarz inequality and letting $k\to \infty$), we have
\[%
\int_{B^{+}_{\rho}} r^{-2}\zeta^{2}\,d\mu_{+} \le C \int_{B^{+}_{\rho}}\bigl(D_{r}\zeta\bigr)^{2}\,d\mu_{+}, %
\quad C=C(\gamma) %
\dl{z-over-r-bd}
\]%
for any $\rho\in [\ha ,1)$ and any Lipschitz $\zeta$ on $B^{+}_{1}$ with support contained in $B^{+}_{\rho}$.

Next observe that $u$ satisfies the weak form of~\ref{equn}:
{\abovedisplayskip8pt\belowdisplayskip8pt%
\[%
\int_{\breve B_{1}^{+}} (u_{r}\zeta_{r}+u_{y}\cdot \zeta_{y})\,d\mu_{+}=0
\dl{weak-form}
\]%
for any Lipschitz $\zeta$ with support in $B^{+}_{\rho}$ for some $\rho<1$, where, here and subsequently,
\[%
d\mu_{+}=r^{\gamma}drdy
\]%
and subscripts denote partial derivatives: 
\[%
u_{r}=D_{r}u,\,\,\,u_{y}= D_{\!y}u=(D_{\!y_{1}}u,\ldots,D_{\!y_{\ell}}u).
\]%
\ref{weak-form} is checked by first observing that it holds with  $\varphi_{\sigma}(r)\zeta(r,y)$ in
place of $\zeta(r,y)$, where $\varphi_{\sigma}(r)=0$ for $r<\sigma/2$, $\varphi_{\sigma}(r)=1$ for $r>\sigma$,
$|\varphi_{\sigma}'(r)|<C/\sigma$, and then letting $\sigma\downarrow 0$ and using~\ref{w-1-2-bd}.

Similarly since $r^{-\gamma}\frac{\partial}{\partial r}(r^{\gamma}\frac{\partial u^{2}}{\partial r})
+\tsum_{j=1}^{\ell}D_{y_{j}}D_{y_{j}}u^{2}= 2(u_{r}^{2}+|u_{y}|^{2})\ge 0$, we can check using~\ref{w-1-2-bd}
and~\ref{z-over-r-bd}  that $u^{2}$ is a weak subsolution, meaning that 
\[%
\int_{\breve B_{1}^{+}}\bigl((u^{2})_{r}\zeta_{r}+(u^{2})_{y}\cdot \zeta_{y}\bigr)\,d\mu_{+}\le 0, 
\dl{subsoln-w}
\]}%
for any non-negative Lipschitz function $\zeta$ on $B^{+}_{1}$ with support $\subset B^{+}_{\rho}$ for some $\rho<1$. 

Next we note that if $\rho\in [\ha,1)$, $r_{0}\ge 0$ and $B_{\rho}^{+}(r_{0},y_{0})\subset B^{+}_{1}$,  then
\[%
|u(r_{0},y_{0})|\le C\bigl(\rho^{-\ell-1-\gamma}\int_{B^{+}_{\rho}(r_{0},y_{0})}u^{2}\,d\mu_{+}\bigr)^{1/2}, %
\quad C=C(\gamma,\ell).
\dl{sup-bd-0}
\]%
To check this, first observe that if  $B^{+}_{\sigma}(r_{0},y_{0})\subset \breve B_{1}^{+}$ and $r_{0}>\sigma$ then the
equation~\ref{equn} is uniformly elliptic divergence form with smooth coefficients on $B_{\sigma/2}(r_{0},y_{0})$, so we can
use standard elliptic estimates for $u$ to give
\[%
|u(r_{0},y_{0})|^{2}\le C\sigma^{-\gamma-1-\ell}\int_{B^{+}_{\sigma}(r_{0},y_{0})}u^{2} \,d\mu_{+}
\dl{sup-bd-0a}
\]%
with $C=C(\gamma,\ell)$. So now assume $B_{\rho}^{+}(r_{0},y_{0})\subset B^{+}_{1}$. If $r_{0}>\rho/4$ we can
use~\ref{sup-bd-0a} with $\sigma=\rho/4< r_{0}$ to give~\ref{sup-bd-0}, while on the other hand if $r_{0}\le \rho/4$ then
we can first take $\sigma=r_{0}/2$ in~\ref{sup-bd-0a} to give
\[%
|u(r_{0},y_{0})|^{2}\le Cr_{0}^{-\ell-1-\gamma}\int_{B^{+}_{r_{0}/2}(r_{0},y_{0})}u^{2}\,d\mu_{+} \le %
Cr_{0}^{-\ell-1-\gamma}\int_{B^{+}_{2r_{0}}(0,y_{0})}u^{2}\,d\mu_{+},  %
\dl{sup-bd-1}
\]%
and then observe, using a straightforward modification of the relevant argument for classical subharmonic functions to the
present case of the $\beta$-subharmonic function $u^{2}$ as in~\ref{subsoln-w}, 
{\abovedisplayskip8pt\belowdisplayskip8pt%
\[%
\text{$\sigma^{-\ell-1-\gamma}\int_{B^{+}_{\sigma}(0,y_{0})}u^{2}\,d\mu_{+}$ is an increasing function  %
of $\sigma$ for $\sigma\in (0,\rho/2]$}.  %
\]%
So from~\ref{sup-bd-1} we conclude
\[%
 |u(r_{0},y_{0})|^{2}\le C \rho^{-\ell-1-\gamma}\int_{B^{+}_{\rho/2}(0,y_{0})}u^{2}\,d\mu_{+} %
\le  C \rho^{-\ell-1-\gamma}\int_{B^{+}_{\rho}(r_{0},y_{0})}u^{2}\,d\mu_{+}
\]%
where $C=C(\gamma,\ell)$.   Thus~\ref{sup-bd-0} is proved.

For $\rho\in [\ha,1)$,  $kh\in \bigl(-(1-\rho),1-\rho\bigr)\setminus\{0\}$, and $k\in \{1,2,\ldots\}$, let $u_{h}^{(k)}$
(defined on $\breve B^{+}_{\rho}$) denote the vector of $k$-th order difference quotients with respect to the $y$-variables;
so for example
\begin{align*}%
u_{h}^{(1)}&=\bigl(h^{-1}(u(x,y+he_{1})-u(x,y)),\ldots,h^{-1}(u(x,y+he_{\ell})-u(x,y))\bigr), \\
u_{h}^{(2)}&=\bigl(h^{-1}(u_{h}^{(1)}(x,y+he_{1})-u_{h}^{(1)}(x,y)),\ldots,
h^{-1}(u_{h}^{(1)}(x,y+he_{\ell})-u_{h}^{(1)}(x,y))\bigr),
\end{align*}}%
 and generally, for $k|h|<1-\rho$,
\begin{align*}%
&u_{h}^{(k)}(x,y) = h^{-1}\bigl((u_{h}^{(k-1)}(x,y+he_{1})-u_{h}^{(k-1)}(x,y)),\ldots, \\
\noalign{\vskip-3pt}
&\hskip2.3in      (u_{h}^{(k-1)}(x,y+he_{\ell})-u_{h}^{(k-1)}(x,y))\bigr), %
\end{align*}%
(which is a function with values in $\R^{\ell^{k}}$).  For notational convenience we also take
\[%
u^{(0)}_{h}(x,y)=u(x,y).
\]%
Then, replacing $\zeta$ in~\ref{weak-form} by $\zeta_{-h}^{(k)}$ and changing variables appropriately (i.e.\ ``integration by
parts'' for finite differences instead of derivatives),
\[%
\int_{\breve B_{1}^{+}}\bigl(D_{r}u^{(k)}_{h}D_{r}\zeta+\tsum_{j=1}^{\ell} D_{y_{j}}u^{(k)}_{h}\, %
D_{y_{j}}\zeta\bigr)\,d\mu_{+}=0,
\]%
for all Lipschitz $\zeta$ on $B_{1}^{+}$ with support $\subset B^{+}_{\rho}$ and for $k|h|<1-\rho$.  Replacing $\zeta$ by
$\zeta^{2}u_{h}^{(k)}$ gives
\begin{align*}%
&\int_{\breve B_{1}^{+}}\bigl(|D_{r}u_{h}^{(k)}|^{2}+ |D_{y}u_{h}^{(k)}|^{2}\bigr)\zeta^{2}  \\  %
\noalign{\vskip-2pt}
&\hskip1in =-2\smash[t]{\int_{\breve B_{1}^{+}}}\bigl(\zeta u_{h}^{(k)}\cdot D_{r}u_{h}^{(k)}D_{r}\zeta + \zeta u_{h}^{(k)}\cdot
D_{r}u_{h}^{(k)}D_{y_{j}}\zeta D_{y_{j}}\zeta\bigr)\,d\mu_{+},%
\end{align*}%
so by Cauchy-Schwarz
\[%
\int_{\breve B_{1}^{+}}\bigl(|D_{r}u_{h}^{(k)}|^{2}+ |D_{y}u_{h}^{(k)}|^{2}\bigr)\zeta^{2}\,d\mu_{+}
\le 4\int_{\breve B_{1}^{+}}|u_{h}^{(k)}|^{2}  |D\zeta|^{2}\,d\mu_{+},%
\]%
and by~\ref{z-over-r-bd} we then have 
\[%
\int_{\breve B_{1}^{+}}\bigl(r^{-2}|u_{h}^{(k)}|^{2}+ |D_{r}u_{h}^{(k)}|^{2}+ %
        |D_{y}u_{h}^{(k)}|^{2}\bigr)\zeta^{2}\,d\mu_{+}  \le C\int_{\breve B_{1}^{+}}|u_{h}^{(k)}|^{2}  %
|D\zeta|^{2}\,d\mu_{+}, 
\dl{u-k-bds-3}
\]%
where $C=C(\gamma,\ell)$.  Now let $u^{(k)}=\lim_{h\to 0}u_{h}^{(k)}$ (i.e.\ $u^{(k)}$ is the array of all mixed partial
derivatives of order $k$ with respect to the $y$ variables defined inductively by $u^{(0)}=u$ and
$u^{(k+1)}=D_{y}u^{(k)}$). Then~\ref{u-k-bds-3} gives
\[%
\int_{\breve B_{1}^{+}}\bigl(r^{-2}|u^{(k)}|^{2}+ |D_{r}u^{(k)}|^{2}+ |u^{(k+1)}|^{2}\bigr)\zeta^{2} %
\le C\int_{\breve B_{1}^{+}}|u^{(k)}|^{2}  |D\zeta|^{2}\,d\mu_{+} %
\dl{u-k-bds-4}
\]%
for each $\zeta\in C^{1}(B^{+}_{\rho})$ with support in $B^{+}_{\theta\rho}$ for some $\theta\in [\ha,1)$ and for each
$k$ such that the right side is finite. By~\ref{w-1-2-bd} the right side is finite for $k=0,1$, and taking $k=1$
in~\ref{u-k-bds-4} then implies the right side is also finite with $k=2$. Proceeding inductively we see that in fact that the
right side is finite for each $k=0,1,\ldots$, so~\ref{u-k-bds-4} is valid and all integrals are finite for all $k=1,2,\ldots$.  

Let $k\in\{1,2,\ldots\}$ and $(r_{0},y_{0})\in B^{+}_{1}$ with $r_{0}\ge 0$.  Then if $B^{+}_{\rho}(r_{0},y_{0})\subset
B^{+}_{1}$ and we let 
\[%
\rho_{j}=\rho-\fr{j}{k}\rho/2,\,\, j=0,\ldots,k-1, 
\]%
and, applying \ref{u-k-bds-4} with $k=j$ and 
\[%
\zeta=1\text{ on }B^{+}_{\rho_{j+1}}(r_{0},y_{0}), \,\, %
\zeta=0\text{ on }\R^{1+\ell}\setminus B^{+}_{\rho_{j}}(r_{0},y_{0}),\text{ and }|D\zeta|\le 3k,  %
\]%
we obtain
\[%
\int_{\breve B^{+}_{\rho_{j+1}}(r_{0},y_{0})} |u^{(j+1)}|^{2}\,d\mu_{+}%
\le C\rho^{-2}k^{2}\int_{\breve B^{+}_{\rho_{j}}(r_{0},y_{0})}|u^{(j)}|^{2}  \,d\mu_{+} %
\]%
with $C=C(\gamma,\ell)$. By iteration this gives
\[%
\int_{\breve B_{\rho/2}^{+}(r_{0},y_{0})} |u^{(k)}|^{2}\,d\mu_{+}%
\le C^{k}\rho^{-2k}(k!)^{2}\int_{\breve B_{\rho}^{+}(r_{0},y_{0})}|u|^{2}  \,d\mu_{+} %
\]%
with suitable $C=C(\gamma,\ell)$ (independent of $k$), and then by~\ref{sup-bd-0} with
$u^{(k)}$ in place of $u$ and $\rho/2$ in place of $\rho$ we obtain
\[%
|u^{(k)}(r_{0},y_{0})|^{2}\le C^{k}\rho^{-2k}(k!)^{2}  %
\rho^{-\gamma-1-\ell}\int_{\breve B^{+}_{\rho}(r_{0},y_{0})}|u|^{2}  \,d\mu_{+}, %
\]%
with $C=C(\gamma,\ell)$ (independent of $k$ and $\rho$).  In view of the arbitrariness of $\rho,r_{0},y_{0}$ this implies 
\[%
\sup_{B^{+}_{\rho/2}} |u^{(k)}|^{2}\le C^{k}(k!)^{2}  %
\int_{\breve B^{+}_{\rho}}|u|^{2}  \,d\mu_{+}, \text{ for each $\rho\in (0,1)$},%
\dl{u-k-bds-5}
\]%
 where $C=C(\gamma,\ell,\rho)$.

Next let $L_{r}$ be the second order operator defined by
\[%
L_{r}(f)=r^{-\gamma}\frac{\partial}{\partial r}\bigl(r^{\gamma}\frac{\partial f}{\partial r}\bigr),
\]%
so that~\ref{equn} says $L_{r}u=-\Delta_{y}u$, where $\Delta_{y}u=\tsum_{j=1}^{\ell}D^{2}_{y_{j}}u$, and by repeatedly
differentiating this identity with respect to the $y$ variables, we have
\[%
L_{r}u^{(k)}(r,y) =-\Delta_{y}u^{(k)},  
\]%
and hence 
\[%
L_{r}^{j}u^{(k)} = -(\Delta_{y})^{j}u^{(k)}. %
\dl{L-j-lap-k}
\]%
for each $j,k=0,1,\ldots$.  In particular, since $|\Delta_{y}^{j}f|\le C^{j}|f^{(2j)}|$ with $C=C(\ell)$, 
\[%
|L_{r}^{j}u^{(k)}|^{2}\le C^{2j} |u^{(k+2j)}|^{2}%
\dl{L-j-k-bd}
\]%
for each $j,k=0,1,2,\ldots$, where $C=C(\ell)$

Next we note that, since $|u^{(k)}|$ is bounded on $B^{+}_{\rho}$ for each $\rho<\ha$ by~\ref{u-k-bds-5}, for small enough
$r>0$ and $|y|<\ha$ we can apply elliptic estimates to give $|D_{r}u^{(k)}| <C/r$ with $C$ fixed independent of $r$ (but
depending on $k$), and since $\gamma>1$ we have $|r^{\gamma}D_{r}u^{(k)}(r,y)|\le Cr^{\gamma-1}\to
0$ as $r\downarrow 0$.  But using~\ref{L-j-k-bd}  with $j=1$ and~\ref{u-k-bds-5} we have
$|D_{r}(r^{\gamma}D_{r}u^{(k)})|\le Cr^{\gamma}$ and hence by integrating with respect to $r$ and using the above fact
that $r^{\gamma}D_{r}u^{(k)}(r,y)\le Cr^{\gamma-1}\to 0$ as $r\downarrow 0$, we then have 
{\abovedisplayskip8pt\belowdisplayskip8pt%
\[%
|D_{r}u^{(k)}(r,y)| \le Cr  \text{ for small enough $r$ and all $|y|\le \rho$, $\rho<\ha$}, 
\dl{u-k-to-0}
\]%
with $C$ depending on $k$, $\gamma$, $\rho$ and $\ell$, and in particular, for each $\rho<\ha$, 
\[%
D_{r}u^{(k)}(0_{+},y)=0=\lim_{r\downarrow 0} D_{r}u^{(k)}(r,y) \text{ uniformly for $|y|\le \rho$.}
\dl{u-k-0-+}
\]%
We now claim the following polynomial approximation property: For each $u$ as in~\ref{equn}, \ref{w-1-2-bd}, each
$j,k=1,2,\ldots$ with $k\le j$,  and each $(r,y)\in\breve B^{+}_{\rho}$, $\rho<\ha$, 
\[%
\bigl|L_{r}^{j-k}u(r,y)-\tsum_{i=1}^{k}c_{ijk}L_{r}^{j-i}u(0,y)r^{2(k-i)}/(2(k-i))!\bigr|\le r^{2k}{\sup}_{B_{\rho}^{+}}
|L_{r}^{j}u|/(2k)!   \dl{inductive-hyp}
\]%
where $0<c_{ijk}\le 1$.  To prove the case $k=1$, by virtue of~\ref{u-k-to-0},\,\ref{u-k-0-+}, we can simply integrate 
from $0$ to $r$, using $D_{r}(r^{\gamma}D_{r}L_{r}^{j-1}u(r,y))=r^{\gamma}L^{j}u(r,y)$, 
followed by a cancellation of $r^{\gamma}$ from each side of the resulting identity. This gives
\[%
|D_{r}L_{r}^{j-1}u(r,y)| \le r\, {\sup}_{B_{\rho}^{+}}|L^{j}u(r,y)|/(\gamma+1),
\]%
and then a second integration using~\ref{u-k-0-+} gives 
\[%
\bigl|L_{r}^{j-1}u(r,y)-L_{r}^{j-1}u(0,y)\bigr|\le (2(\gamma+1))^{-1}r^{2}
{\sup}_{B_{\rho}^{+}}|L_{r}^{j}u(r,y)|,\,\,j=1,2,\ldots,%
\]%
which establishes the case $k=1$.  Assume $k+1\le j$ and that~\ref{inductive-hyp} is correct for $k$.  Multiplying each side
of~\ref{inductive-hyp} by $r^{\gamma}$ and integrating, we obtain
\begin{align*}%
&\bigl|r^{\gamma}D_{r} L_{r}^{j-k-1}u(r,y)  %
 -\tsum_{i=1}^{k}(2(k-i)+\gamma+1)^{-1}    c_{ijk}L_{r}^{j-i}u(0,y)\frac{r^{2(k-i)+\gamma+1}}{(2(k-i))!}\bigr| \\  %
\noalign{\vskip-1pt} %
&\hskip1.8in \le (2k+\gamma+1)^{-1}\frac{r^{2k+\gamma+1}}{(2k)!}{\sup}_{B_{\rho}^{+}}|L_{r}^{j}u|,  
\end{align*}%
where we used the fact that $r^{\gamma}L_{r}^{j-k}u(r,y)=D_{r}(r^{\gamma}D_{r}L_{r}^{j-k-1}u(r,y))$.

\vskip1pt

After cancelling the factor $r^{\gamma}$ and integrating again, we obtain
\begin{align*}%
&\bigl| L_{r}^{j-k-1}u(r,y)-L_{r}^{j-k-1}u(0,y)  \\  %
&\hskip0.1in -\tsum_{i=1}^{k}(2(k-i)+\gamma+1)^{-1}(2(k-i)+2)^{-1}   %
                                                 c_{ijk}L_{r}^{j-i}u(0,y)\frac{r^{2(k-i)+2}}{(2(k-i))!}\bigr|\\ %
&\hskip0.15in \le (2k+\gamma+1)^{-1}(2k+2)^{-1}r^{2k+2}  %
                     {\sup}_{B_{\rho}^{+}} |L_{r}^{j}u|/(2k)!   \le  {\sup}_{B_{\rho}^{+}} |L_{r}^{j}u|r^{2k+2}/(2k+2)! %
\end{align*}%
which confirms the validity of~\ref{inductive-hyp} with $k+1$ in place of $k$. So~\ref{inductive-hyp} is proved for all $k\le
j$, and in particular with $k=j$ and suitable constants $c_{ij}\in (0,1]$ we get
\[%
u(r,y)=\tsum_{i=0}^{j-1}c_{ij}L_{r}^{i}u(0,y)r^{2i}/(2i)!+E_{j}(r,y),\text{where\,} |E_{j}(r,y)|\le %
                                                                                  r^{2j}{\sup}_{B_{\rho}^{+}} |L_{r}^{j}u|/(2j)!   %
\]%
By~\ref{L-j-k-bd} and~\ref{u-k-bds-5}, 
\[%
{\sup}_{B_{\rho/2}^{+}}|L_{r}^{i}u|/(2i)!={\sup}_{B_{\rho/2}^{+}}|\Delta_{y}^{i}u|/(2i)!\le  %
C^{i}{\sup}_{B_{\rho/2}^{+}}|u^{(2i)}|/(2i)!\le C^{i} \bigl(\int_{\breve B^{+}_{\rho}}|u|^{2} \,d\mu_{+}\bigr)^{1/2},   %
\]%
with $C=C(\gamma,\ell,\rho)$, so, for suitable $C=C(\gamma,\ell)$, we conclude that $u(r,y)$ has a power series
expansion in terms $r^{2}$:
\[%
u(r,y) =\tsum_{j=0}^{\infty}\,\,\,\,a_{j}(y)r^{2j}, \quad 0\le r< \sigma,  
\dl{r-exp}
\]%
where $\sigma=\sigma(\gamma,\ell)\in (0,\ha]$, and $a_{j}$ satisfies the bounds
\[%
{\sup}_{B^{+}_{\sigma}}|a_{j}|\le
C^{j}\bigl(\int_{B_{1/2}^{+}}u^{2}\,d\mu_{+}\bigr)^{1/2}, 
\]%
where $C=C(\gamma,\ell)$.  Thus~\ref{r-exp} implies
\[%
{\sup}_{B_{\sigma/2}^{+}}|D_{r}^{j}u(r,y)| \le C^{j}j! \,\bigl(\int_{B_{1/2}^{+}}u^{2}\,d\mu_{+}\bigr)^{1/2}, %
  \,\,\,C=C(\gamma,\ell). %
\dl{bds-D-j}
\]%
Since the same holds with $u^{(k)}$ in place of $u$,  and since
\[%
\int_{B_{1/2}^{+}}|u^{(k)}|^{2}\,d\mu_{+}\le C^{k}(k!)^{2}\int_{B_{3/4}^{+}}u^{2}\,d\mu_{+}, 
\]}%
by~\ref{u-k-bds-5}, we deduce from~\ref{bds-D-j} that for suitable $\sigma=\sigma(\gamma,\ell)\in (0,\ha)$
\[%
{\sup}_{B_{\sigma}^{+}}
|D_{r}^{j}D_{y}^{k}u(r,y)| \le C^{j+k} j!k! \,\bigl(\int_{B_{3/4}^{+}}u^{2}\,d\mu_{+}\bigr)^{1/2},  %
\dl{deriv-bds}
\]%
where $C=C(\gamma,\ell)$, and hence in particular $u$ is real-analytic in the variables $r^{2}$ and $y_{1},\ldots,y_{\ell}$
in a neighborhood of $(0,0)$ as claimed.

Finally we show that if $u$ is $\beta$-harmonic in $\breve B_{1}^{+}$ then the power series for $u$ converges in
$B_{\rho}^{+}$ for each $\rho<1$, and also that the homogeneous $\beta$-harmonic polynomials restricted to $S_{+}^{\ell}$
are complete in $L^{2}(\nu_{+})$ on $S_{+}^{\ell}$, where $\nu_{+}$ is the measure
$d\nu_{+}=\omega_{1}^{\gamma}d\mu_{\ell}$ on $S_{+}^{\ell}$.

So let $u\in L^{2}(\mu_{+})$ satisfy~\ref{equn} and \ref{w-1-2-bd}.  The above discussion shows that for suitably small
$\sigma$ we can write
\[%
u=\tsum_{j=0}^{\infty}u_{j} \,\,\,\text{ in }B_{\sigma}^{+}, %
\dl{rau}
\]%
where $u_{j}$ consists of the homogeneous degree $j$ terms in the power series expansion of $u$ in $B_{\sigma}$ (and
$u_{j}=0$ if there are no such terms). Then each $u_{j}\neq 0$ is a homogeneous degree j $\beta$-harmonic polynomial and
we let
$$%
\tilde u_{j}(\omega)=\rho^{-j}u_{j}(\rho\omega),\,\,\, \hat u_{j}(\omega)=
\|\tilde u_{j}\|_{L^{2}(\nu_{+})}^{-1}\tilde u_{j}(\omega),   \,\,\, \omega\in \Sph^{\ell}_{+},
$$%
and we set $\hat u_{j}(\omega)=0$ if $u_{j}=0$.  Then, with $\langle\,,\,\rangle=$ the $L^{2}(\nu_{+})$ inner product,
$\langle u,\hat u_{j} \rangle \hat u_{j}=u_{j}$ in $B^{+}_{\sigma}$ for each $j$, and hence by~\ref{rau} the series
$\sum_{j}\langle u,\hat u_{j} \rangle \hat u_{j}$ converges smoothly (and also in $L^{2}(\mu_{+})$) to $u$ in
$B^{+}_{\sigma}$. 

By definition $\rho^{j}\hat u_{j}(\omega)$ is either zero or a homogeneous degree $j$ harmonic polynomial, so
by~\ref{jth-eig} \smash{$\dvg_{\Sph^{\ell}_{+}}(\omega_{1}^{\gamma}\nabla_{\Sph^{\ell}_{+}}\hat u_{j})=- j(j+\ell+\beta)
\omega_{1}^{\gamma}\hat u_{j}$}, and hence using the formula~\ref{sph-coords} we can directly check that $\langle u,\hat
u_{j} \rangle \hat u_{j}$ is $\beta$-harmonic on all of $\breve B_{1}^{+}$. Since by construction it is equal to $u_{j}$ on
$B_{\sigma}$, by unique continuation (applicable since $u$ is real-analytic on $\breve B_{1}^{+}$) we conclude
\begin{align*}%
&\langle u,\hat u_{j} \rangle \hat u_{j} \text{ is either zero or the homogeneous degree $j$ $\beta$-harmonic} \dtg{poly} \\  %
\noalign{\vskip-5pt}
&\hskip1.4in\text{ polynomial $u_{j}$ on all of $B_{1}^{+}\setminus S^{\ell}_{+}$ for each $j=0,1,\ldots$}.
\end{align*}%
Also, by the orthogonality~\ref{orthog},
\begin{align*}%
\bigl\|\sum_{j=p}^{q}\bigl\langle u,\hat u_{j} \bigr\rangle \hat u_{j}\bigr\|^{2}_{L^{2}(\mu_{+}^{\rho})} %
&=\sum_{j=p}^{q}\int_{0}^{\rho}\bigl\langle u(\tau\omega),\hat u_{j}(\omega)\bigr\rangle^{2}\, %
\tau^{\gamma+\ell}d\tau \\   %
&\le \int_{0}^{\rho}\bigl\|u(\tau\omega)\bigr\|^{2}_{L^{2}(\nu_{+})}\,\tau^{\gamma+\ell}d\rho=  %
\bigl\|u\bigr\|^{2}_{L^{2}(\mu^{\rho}_{+})} \,\,(<\,\infty)  %
\end{align*}%
for each $\rho<1$ and each $p<q$, where $\smash{\mu_{+}^{\rho}}$ is the measure $\mu_{+}$ on $B_{\rho}$.  So
\smash{$\sum_{j=0}^{q}\langle u,\hat u_{j} \rangle \hat u_{j}$} is Cauchy, hence convergent, in \smash{$L^{2}(\mu_{+}^{\rho})$}
to a $\beta$-harmonic function $v$ on $\breve B_{\rho}^{+}$.  But $v=u$ on $B_{\sigma}^{+}$ and hence, again using unique
continuation, $v=u$ in all of $\breve B_{\rho}^{+}$.  Thus \smash{$\sum_{j=0}^{q}\langle u,\hat u_{j} \rangle \hat u_{j}$}
converges to $u$ in $L^{2}(\smash{\mu_{+}^{\rho}})$ for each $\rho<1$ and the convergence is in $L^{2}(\mu_{+})$ if
$\|u\|_{L^{2}(\mu_{+})}<\infty$.

\vskip1pt

Now observe that the bounds~\ref{deriv-bds} were established for balls centred at $(0,0)$, but with only notational changes
the same argument gives similar bounds in balls centred at $(0,y_{0})$ with $|y_{0}|<1$. Specifically for each
$\rho\in (0,1)$ and each 
$|y_{0}|<\rho$ there is $\sigma=\sigma(\gamma,\ell,\rho)<\ha(1-\rho)$ such that 
\[%
{\sup}_{B_{\sigma}^{+}(0,y_{0})}|D_{r}^{j}D_{y}^{k}u|\le
C^{j+k}j!k!\bigl(\int_{B_{(1-\rho)/2}^{+}(0,y_{0})}u^{2}\,d\mu_{+}\bigr)^{1/2},\,\,C=C(\gamma,\ell,\rho).
\]%
So in fact,  with $\sigma=\sigma(\gamma,\ell,\rho)$ small enough,
\[%
{\sup}_{\{(r,y):r\in [0,\sigma],|y|\le \rho\}}
|D_{r}^{j}D_{y}^{k}u(r,y)| \le C^{j+k} j!k! \bigl(\int_{B_{1}^{+}}u^{2}\,d\mu_{+}\bigr)^{1/2},\quad \rho<1.  %
\]%
Also in $B^{+}_{\rho}\setminus ([0,\sigma]\times \R^{\ell})$ we can use standard elliptic estimates,  so in fact we have
\[%
{\sup}_{B^{+}_{\rho}}
|D_{r}^{j}D_{y}^{k}u(r,y)| \le C \bigl(\int_{B_{1}^{+}}u^{2}\,d\mu_{+}\bigr)^{1/2}, %
\dl{mod-deriv-bds-2}
\]%
with $C=C(j,k,\gamma,\rho,\ell)$, so the $L^{2}$ convergence of the series $\sum_{j}\langle u,\hat u_{j} \rangle \hat
u_{j}(=\sum_{j}u_{j})$ proved above is also $C^{k}$ convergence in $B_{\rho}^{+}$ for each $k\ge 1$ and each $\rho<1$. 

\vskip1pt

Finally to prove the completeness of the homogeneous $\beta$-harmonic polynomials in $L^{2}(\nu_{+})$ (on
$\Sph^{\ell}_{+}$), let $\varphi$ be any smooth function on $\Sph^{\ell}_{+}$ with $\varphi$ zero in some neighborhood of
$r=0$. By minimizing the energy \smash{$\int_{B_{1}^{+}}(u_{r}^{2}+|u_{y}|^{2})\,r^{\gamma}d\mu$} among functions with trace
$\varphi$ on $\Sph_{\ell}$ we obtain a solution of~\ref{weak-form} with trace $\varphi$ on $\Sph^{\ell}_{+}$.  The above
discussion plus elliptic boundary regularity shows that $u$ is $C^{0}$ on all of $B_{1}^{+}$ and that the sequence
\smash{$\{\sum_{j=0}^{q}u_{j}\}_{q=0,1,2,\ldots}$}, which we showed above to be convergent to $u$ in $L^{2}(\mu_{+})$ on
$B_{1}^{+}$, is also uniformly convergent to $u$ on all of $B_{1}^{+}$. Hence
\smash{$\varphi(\omega)=\sum_{j=0}^{\infty}u_{j}(\omega)$} on $\Sph^{\ell}_{+}$ with the convergence uniform and hence in
$L^{2}(\nu_{+})$.  Thus $\varphi$ is represented as an $L^{2}(\nu_{+})$ convergent series of $\beta$-harmonic polynomials
on $\Sph^{\ell}_{+}$.  Since such $\varphi$ are dense in $L^{2}(\nu_{+})$, the required completeness is established.


\bibliographystyle{amsalpha}


\newcommand{\noopsort}[1]{}
\providecommand{\bysame}{\leavevmode\hbox to3em{\hrulefill}\thinspace}
\providecommand{\MR}{\relax\ifhmode\unskip\space\fi MR }
\providecommand{\MRhref}[2]{%
  \href{http://www.ams.org/mathscinet-getitem?mr=#1}{#2}
}
\providecommand{\href}[2]{#2}

\end{document}